%
%
%
%
\hsize=5in
\baselineskip=12pt
\vsize=20cm
\parindent=10pt
\pretolerance=40
\predisplaypenalty=0
\displaywidowpenalty=0
\finalhyphendemerits=0
\hfuzz=2pt
\frenchspacing
\footline={\ifnum\pageno=1\else\hfil\tenrm\number\pageno\hfil\fi}
%
%
\input amssym.def
\def\titlefonts{\baselineskip=1.44\baselineskip
	\font\titlef=cmbx12
	\titlef
	}
\font\tenib=cmmib10 
\font\sevenib=cmmib7
\skewchar\tenib='177
\skewchar\sevenib='177
\font\tenbsy=cmbsy10 
\font\sevenbsy=cmbsy7
\skewchar\tenbsy='60
\skewchar\sevenbsy='60
\def\boldfonts{\bf
	\textfont1=\tenib
	\scriptfont1=\sevenib
	\textfont2=\tenbsy
	\scriptfont2=\sevenbsy
	}
\font\ninerm=cmr9
\font\ninebf=cmbx9
\font\ninei=cmmi9
\skewchar\ninei='177
\font\ninesy=cmsy9
\skewchar\ninesy='60
\font\nineit=cmti9
\def\reffonts{\baselineskip=0.9\baselineskip
	\textfont0=\ninerm
	\def\rm{\fam0\ninerm}%
	\textfont1=\ninei
	\textfont2=\ninesy
	\textfont\bffam=\ninebf
	\def\bf{\fam\bffam\ninebf}%
	\def\it{\nineit}%
	}
%
%
\def\frontmatter{\vbox{}\vskip1cm\bgroup
	\leftskip=0pt plus1fil\rightskip=0pt plus1fil
	\parindent=0pt
	\parfillskip=0pt
	\pretolerance=10000
	}
\def\endfrontmatter{\egroup\bigskip}
\def\title#1{{\titlefonts#1\par}}
\def\author#1{\bigskip#1\par}
\def\address#1{\medskip{\reffonts\it#1}}
\def\email#1{\smallskip{\reffonts{\it E-mail: }\rm#1}}
\def\thanks#1{\footnote{}{\reffonts\rm\noindent#1\hfil}}
\def\section#1\par{\ifdim\lastskip<\bigskipamount\removelastskip\fi
	\penalty-250\bigskip
	\vbox{\leftskip=0pt plus1fil\rightskip=0pt plus1fil
	\parindent=0pt
	\parfillskip=0pt
	\pretolerance=10000{\boldfonts#1}}\nobreak\medskip
	}
\def\emph#1{{\it#1}\/}
\def\proclaim#1. {\medbreak\bgroup{\noindent\bf#1.}\ \it}
\def\endproclaim{\egroup
	\ifdim\lastskip<\medskipamount\removelastskip\medskip\fi}
\newdimen\itemsize
\def\setitemsize#1 {{\setbox0\hbox{#1\ }
	\global\itemsize=\wd0}}
\def\item#1 #2\par{\ifdim\lastskip<\smallskipamount\removelastskip\smallskip\fi
	{\leftskip=\itemsize
	\noindent\hskip-\leftskip
	\hbox to\leftskip{\hfil\rm#1\ }#2\par}\smallskip}
\def\Proof#1. {\ifdim\lastskip<\medskipamount\removelastskip\medskip\fi
	{\noindent\it Proof\if\space#1\space\else\ \fi#1.}\ }
\def\endproof{\hfill\hbox{}\quad\hbox{}\hfill\llap{$\square$}\medskip}
\def\Remark. {\ifdim\lastskip<\medskipamount\removelastskip\medskip\fi
        {\noindent\bf Remark. }}
\def\endremark{\medskip}
%
%
\newcount\citation
\newtoks\citetoks
\def\citedef#1\endcitedef{\citetoks={#1\endcitedef}}
\def\endcitedef#1\endcitedef{}
\def\citenum#1{\citation=0\def\curcite{#1}%
	\expandafter\checkendcite\the\citetoks}
\def\checkendcite#1{\ifx\endcitedef#1?\else
	\expandafter\lookcite\expandafter#1\fi}
\def\lookcite#1 {\advance\citation by1\def\auxcite{#1}%
	\ifx\auxcite\curcite\the\citation\expandafter\endcitedef\else
	\expandafter\checkendcite\fi}
\def\cite#1{\makecite#1,\cite}
\def\makecite#1,#2{[\citenum{#1}\ifx\cite#2]\else\expandafter\clearcite\expandafter#2\fi}
\def\clearcite#1,\cite{, #1]}
%
%
\def\references{\section References\par
	\bgroup
	\parindent=0pt
	\reffonts
	\rm
	\frenchspacing
	\setbox0\hbox{99. }\leftskip=\wd0
	}
\def\endreferences{\egroup}
\newtoks\authtoks
\newif\iffirstauth
\def\checkendauth#1{\ifx\auth#1%
    \iffirstauth\the\authtoks
    \else{} and \the\authtoks\fi,%
  \else\iffirstauth\the\authtoks\firstauthfalse
    \else, \the\authtoks\fi
    \expandafter\nextauth\expandafter#1\fi
}
\def\nextauth#1,#2;{\authtoks={#1 #2}\checkendauth}
\def\auth#1{\nextauth#1;\auth}
\newif\ifbookref
\def\nextref#1 {\par\hskip-\leftskip
	\hbox to\leftskip{\hfil\citenum{#1}.\ }%
	\initnextref}
\def\initnextref{\bookreffalse\firstauthtrue\ignorespaces}
\def\paper#1{{\it#1},}
\def\book#1{\bookreftrue{\it#1},}
\def\journal#1{#1}
\def\BkSer#1{#1,}
\def\Vol#1{{\bf#1}}
\def\BkVol#1{Vol. #1,}
\def\publisher#1{#1,}
\def\Year#1{\ifbookref #1.\else(#1)\fi}
\def\Pages#1{\makepages#1.}
\long\def\makepages#1-#2.#3{#1--#2\ifx\par#3.\else; \fi}
\def\inRus{{ \rm(in Russian)}}
\def\etransl#1{English translation in \journal{#1}}
%
%
\newsymbol\varnothing 203F
\newsymbol\diagdown 231F
\newsymbol\diagup 231E
\newsymbol\square 1203
\let\ot\otimes
\let\sbs\subset
\let\<\langle
\let\>\rangle
\def\alt{_{\rm alt}}
\def\ch{\mathop{\rm ch}\nolimits}
\def\chr{\mathop{\rm char}\nolimits}
\def\Coend{\mathop{\rm Coend}\nolimits}
\def\deg{\mathop{\rm deg}\nolimits}
\def\dimk{\dim\nolimits_{\mskip2mu\bbk}}
\def\End{\mathop{\rm End}\nolimits}
\def\frac#1#2{{#1\over#2}}
\def\Grot{\mathop{\rm Grot}\nolimits}
\def\Grotr{\mathop{\rm Grot}\nolimits^r}
\def\Hom{\mathop{\rm Hom}\nolimits}
\def\Homk{\Hom\nolimits_{\mskip2mu\bbk}}
\def\HomZ{\Hom\nolimits_{\mskip2mu\bbZ}}
\def\Id{\mathop{\rm Id}\nolimits}
\def\Im{\mathop{\rm Im}}
\def\ind{\mathop{\rm ind}\nolimits}
\def\Ker{\mathop{\rm Ker}}
\def\mod#1{\ifinner\mskip8mu(\mathop{\rm mod}#1)
        \else\mskip12mu(\mathop{\rm mod}#1)\fi}
\def\Sym{\mathop{\rm Sym}\nolimits}
\def\tprod{\textstyle\prod\limits}
\def\triv{_{\rm triv}}
\def\tsum{\textstyle\sum}
\def\mapr#1{{}\mathrel{\smash{\mathop{\longrightarrow}\limits^{#1}}}{}}
\def\lmapr#1#2{{}\mathrel{\smash{\mathop{\count0=#1
  \loop
    \ifnum\count0>0
    \advance\count0 by-1\smash{\mathord-}\mkern-4mu
  \repeat
  \mathord\rightarrow}\limits^{#2}}}{}}
\def\mapse#1{\vcenter{\offinterlineskip
  \setbox0\hbox{$\diagdown$}\dimen0=\wd0
  \vss\box0
  \vskip-0.1\dimen0
  \hbox{\hskip0.82\dimen0
    \raise0.2\dimen0\hbox{\llap{$\scriptstyle#1$\hskip0.3\dimen0}}$\searrow$}
  \vss}
}
\def\mapne#1{\vcenter{\offinterlineskip
  \setbox0\hbox{$\diagup$}\dimen0=\wd0
  \vss\hbox{\hskip0.80\dimen0$\nearrow$}
  \vskip-0.1\dimen0
    \hbox{\box0\raise0.2\dimen0\hbox{\rlap{$\hskip0.4\dimen0\scriptstyle#1$}}}
  \vss}
}
\def\diagram#1{\vbox{\halign{&\hfil$##$\hfil\cr #1}}}
\let\al\alpha
\let\be\beta
\let\ga\gamma
\let\de\delta
\let\la\lambda
\let\om\omega
\let\ph\varphi
\let\si\sigma
\let\th\theta
\let\ze\zeta
\let\Ga\Gamma
\let\De\Delta
\let\La\Lambda
\let\Si\Sigma
\def\frm{{\frak m}}
\def\frA{{\frak A}}
\def\frS{{\frak S}}
\def\calD{{\cal D}}
\def\calH{{\cal H}}
\def\calI{{\cal I}}
\def\calR{{\cal R}}
\newsymbol\bbk 207C
\def\bbC{{\Bbb C}}
\def\bbN{{\Bbb N}}
\def\bbR{{\Bbb R}}
\def\bbS{{\Bbb S}}
\def\bbT{{\Bbb T}}
\def\bbZ{{\Bbb Z}}
\def\fhat{{\hat f}}
\def\ghat{{\hat g}}
\def\phat{{\hat p}}
\def\Lov{{\overline L}}

\def\Ak{A_{\mskip0.5mu\bbk}}
\def\chilamu{\chi\strut_\lamu}

\def\chR{\ch^{\mskip-1mu R}}
\def\dla{d^{\mskip1mu\la}}
\def\drho{d^{\mskip1mu\rho}}
\def\EndHM{\End_{\calH_n(z)}\mskip-2mu M}
\def\EndHMk{\End_{\calH_n(q)}\mskip-2mu M_\bbk}
\def\EndHMQ{\End_{\calH_n(z)_Q}\mskip-2mu M_Q}
\def\EndHn{\End_{\calH_n}\mskip-5mu}
\def\EndHX{\End_{\calH_n(q)}\mskip-2mu X}
\def\Endk{\End\nolimits_{\mskip2mu\bbk}}
\def\EndkV{\Endk\mskip-2mu V}
\def\iin{_{i,\mskip1mu n-i}}
\def\lai{{\la,\mskip1mu i}}
\def\lamu{{\la,\mskip1mu\mu}}
\def\lazero{{\la,\mskip1mu0}}
\def\Lk{L_{\mskip1mu\bbk}}
\def\Mk{M_{\mskip1mu\bbk}}
\def\Nk{N_{\mskip1mu\bbk}}
\def\Nlamunu{N_{(\lamu),\mskip1mu\nu}}
\def\Par{{\cal P}}
\def\psik{\psi_{\mskip1mu\bbk}}
\def\Rep{\frak{Rep}^1}
\def\RepA{\mathop{\rm Rep}\frA}
\def\RepHq{{\rm Rep}^1\mskip1mu\calH_n(q)}
\def\RHq{\mathop{\rm Rep}\calH_n(q)}
\def\rhth{{\rho,\mskip1mu\th}}
\def\Smu{\bbS\vbox to\ht\strutbox{}^{\mskip0.5mu\mu}}
\def\Smula{\Smu_{\mskip1mu\la}}
\def\Triv{\frak{Triv}}
\def\TrivHq{\mathop{\rm Triv}\calH_n(q)}
\def\UK{U_{\mskip-2mu K}}
\def\zeromu{{0,\mskip1mu\mu}}
\def\PHHai{Ph\`ung H\^o Hai }

\citedef
Ais-SW52
At-T69
Ber-R87
Br-K03
Dav00
Dip-J86
Dip-J89
Dip-J91
Don
Don01
Du-PW91
Edr52
Geck-P
Gia-LDW17
Gur90
Hai99
Hai02
Hai05
Hen
Lar-T91
Mac
Resh-TF89
Ser84
Sk19
Stan
Zel
\endcitedef

\frontmatter

\title{On the graded algebras associated with Hecke symmetries, II. The Hilbert series}
\author{Serge Skryabin}
\address{Institute of Mathematics and Mechanics,
Kazan Federal University,\break
Kremlevskaya St.~18, 420008 Kazan, Russia}
\email{Serge.Skryabin@kpfu.ru}

\endfrontmatter

\section
Introduction

With a Hecke symmetry $R$ on a finite dimensional vector space $V$ one 
associates the $R$-symmetric algebra $\bbS(V,R)$, the $R$-skewsymmetric 
algebra $\La(V,R)$, and the bialgebra $A(R)$ given by the 
Faddeev-Reshetikhin-Takhtajan construction \cite{Resh-TF89}. The first two 
algebras are noncommutative analogs of the symmetric and the exterior algebras 
of $V$. The algebra $A(R)$ is a noncommutative analog of the ring of 
polynomial functions on the space of $r\times r$ matrices. By localizing 
$A(R)$ one obtains a Hopf algebra which represents a nonstandard quantum group 
(see \cite{Gur90}).

For two Hecke symmetries $R,R'$ with the same parameter $q$ of the Hecke 
relation there is also an algebra $A(R',R)$ generalizing $A(R)$. This algebra 
introduced by \PHHai \cite{Hai02} in a different notation represents a 
``quantum hom-space". All these algebras are quadratic graded algebras.

The present paper continues the work started in \cite{Sk19} which attempts to 
understand the general properties of the graded algebras associated with 
Hecke symmetries without a restriction on the parameter $q$ of the Hecke 
relation imposed in earlier results. The known results rely heavily on 
semisimplicity of the Hecke algebras $\calH_n=\calH_n(q)$ of type $A$ which 
operate in the tensor powers of the initial space $V$. This \emph{semisimple 
case} occurs precisely when
$$
1+q+\ldots+q^{n-1}\ne0\quad\hbox{for all $\,n>0$}.
$$

Here we will be concerned with the determination of the Hilbert series and 
several related results. For a graded algebra $A=A_0\oplus A_1\oplus 
A_2\oplus\ldots$ with finite dimensional homogeneous components its 
\emph{Hilbert series} is a formal power series in one indeterminate defined as
$$
H_A(t)=\sum\,(\dim A_n)\,t^n.
$$

The question as to what are possible Hilbert series of the algebras $\bbS(V,R)$ 
and $\La(V,R)$ arose in the work of Gurevich \cite{Gur90}. About 10 years 
later \PHHai \cite{Hai99} and, independently, Davydov \cite{Dav00} observed 
that in the semisimple case the dimensions of the homogeneous components of 
the two algebras form totally positive sequences. From this they deduced that 
the Hilbert series of these algebras are rational functions with negative 
roots and positive poles. This conclusion is based on analytic results 
obtained by Aissen, Schoenberg, Whitney \cite{Ais-SW52} and Edrei \cite{Edr52} 
which describe the generating series of totally positive sequences.

As we have seen in \cite{Sk19}, good properties may be lost when $q$ is a root 
of 1. However, it was shown there that several previously known results extend 
to the case of an arbitrary $q$ provided that a certain additional condition 
is imposed. Recall that an indecomposable $\calH_n$-module is said to have a 
\emph{$1$-dimensional source} if it is a direct summand of an $\calH_n$-module 
induced from a $1$-dimensional representation of a parabolic subalgebra, and 
we say that $R$ satisfies the \emph{$1$-dimensional source condition} if for 
each $n>0$ all indecomposable direct summands of $\,V^{\ot n}$ regarded as an 
$\calH_n$-module with respect to the representation arising from $R$ have 
$1$-dimensional sources. This condition is satisfied automatically in the 
semisimple case mentioned earlier. Our main result in the new paper is

\proclaim
Theorem 3.8.
Suppose that $R$ satisfies the $1$-dimensional source condition. Then
$$
H_{\La(V,R)}(t)=f_0(-t)/f_1(t),\qquad H_{\bbS(V,R)}(t)=f_1(-t)/f_0(t)
$$
with integer polynomials $\,f_0,\,f_1\in\bbZ[t]\,$ whose constant terms are 
equal to\/ $1$ and all roots are positive real numbers.
\endproclaim

The pair $(r_0,r_1)$ where $r_i=\deg f_i$ for $i=0,1$ is called the 
\emph{birank} of $R$. Thus the Hilbert series of the two algebras can be 
written as
$$
{\tprod_{i=1}^{r_0}}(1+\al_it)\cdot
{\tprod_{j=1}^{r_1}}(1-\be_jt)^{-1}\quad{\rm and}\qquad
{\tprod_{j=1}^{r_1}}(1+\be_jt)\cdot
{\tprod_{i=1}^{r_0}}(1-\al_it)^{-1}\quad
$$
where $\al_i$ and $\be_j$ are positive algebraic integers.


The already mentioned results of \PHHai and Davydov may be viewed as a nice 
application of the theory of symmetric functions. We will use a nonstandard 
notation $\,\Sym\,$ for the ring of symmetric functions defined as in Macdonald 
\cite{Mac} (in our paper the letter $\La$ is reserved for a different ring). 
Consider the Grothendieck ring $\Grot(R)$ of the category of finite dimensional 
right $A(R)$-comodules. In the semisimple case one can use a quantum version of 
the Schur-Weyl duality to obtain a ring homomorphism $\Sym\to\Grot(R)$ under 
which each Schur function $s_\la$ is sent either to 0 or to the class of a 
simple comodule.

Since $\Sym$ is a polynomial ring in a countable set of indeterminates, a 
homomorphism $\ph:\Sym\to\Grot(R)$ can be easily constructed in the case of 
arbitrary $q$ by specifying its values on the generators. The main obstacle we 
encounter is to show that the element $\ph(s_\la)\in\Grot(R)$ is \emph{positive} 
in the sense that $\ph(s_\la)$ is the class of an actual comodule $V^\la$, in 
general not defined uniquely, of course. What is needed here can be reformulated 
in terms of the Grothendieck group of the category of finite dimensional 
$\EndHn X$-modules where $X=V^{\ot n}$ with the $\calH_n$-module structure 
arising from $R$. The necessary property is stated in Corollary 2.8, and 
section 2 is devoted to its proof. A key role is played by a version of the 
decomposition map which provides a bridge between the Grothendieck groups in
the semisimple and nonsemisimple cases.

Total positivity of the sequence $\bigl(\dim\bbS_n(V,R)\bigr)$ is an immediate 
consequence of positivity of the images of the Schur functions under $\ph$. 
Once it is known, we can invoke the analytic result of \cite{Ais-SW52} and 
\cite{Edr52}. Actually it will be shown in section 3 that rationality of the 
Hilbert series can be explained by purely algebraic arguments, and the 
remainder of the proof is then much shorter than in general. In this way we 
present a selfcontained proof of Theorem 3.8.

Under the same assumption about $R$ it will be shown in section 4 that the 
class of $V^{\ot n}$ in the Grothendieck group of the category 
of finite dimensional $\calH_n$-modules is completely determined by the 
Hilbert series of $\bbS(V,R)$. Moreover, we describe in Theorem 4.5 a certain 
element $\ch(V^{\ot n})\in\Sym$ which contains full information about this 
class $[V^{\ot n}]$. However, it is not clear whether $V^{\ot n}$ can always 
be determined as an $\calH_n$-module up to isomorphism.

If the algebra $\La(V,R)$ is finite dimensional and $R$ satisfies the 
\emph{trivial source condition} in the sense that for each $n>0$ the 
indecomposable $\calH_n$-module direct summands of $V^{\ot n}$ are induced 
from the \emph{trivial} $1$-dimensional representations of parabolic 
subalgebras of $\calH_n$, then our results are much more complete. Indeed, the 
$\calH_n$-module $V^{\ot n}$ is described in Theorem 6.1. As a consequence, in 
this case the algebra $A_n(R)^*$ dual to the subcoalgebra $A_n(R)\sbs A(R)$ is 
Morita equivalent to the $q$-Schur algebra of Dipper and James $S_q(r,n)$ 
\cite{Dip-J91} where $r=r_0$ is the \emph{rank} of $R$ (we denote by $A_n(R)$ 
the degree $n$ homogeneous component of $A(R)$). Therefore the category of 
$A_n(R)$-comodules is equivalent to the well studied category of 
$S_q(r,n)$-modules. In particular, this category depends only on $q$, $r$, and 
$n$, but not on $R$ itself. It is a highest weight category (see Donkin 
\cite{Don}).

It has been known for a long time that $A_n(R)^*\cong S_q(r,n)$ for many 
different quantizations of the semigroup of $r\times r$ matrices. This 
phenomenon was first observed by Du, Parshall and Wang \cite{Du-PW91} in the 
case of Takeuchi's 2-parameter family of deformations. In an equivalent 
formulation, two bialgebras $A(R)$ and $A(R')$ in this family are isomorphic 
as coalgebras whenever the corresponding parameters satisfy a certain 
relation, and then their corepresentation categories are obviously equivalent. 
As is seen from \cite{Du-PW91, (2.7)} it was not clear at that time whether 
these two categories are monoidally equivalent. On the level of Hopf envelopes 
a general result on braided monoidal equivalence was obtained later by \PHHai 
in the semisimple case \cite{Hai05}.

We will use Theorem 6.1 to strengthen two results from the previous paper 
\cite{Sk19}. Keeping the previous assumption about $R$, let $R'$ be a second 
Hecke symmetry satisfying the same conditions. Theorem 6.3 states that there 
is a braided monoidal equivalence between the categories of $A(R)$-comodules 
and $A(R')$-comodules provided that the two Hecke symmetries have the same 
parameter $q$ and the same rank $r$. This equivalence is obtained by cotensoring 
right $A(R)$-comodules with the bicomodule algebra $A(R',R)$ (see \cite{Sk19, 
Th. 7.2}). By Theorem 6.4 the graded algebra $A(R',R)$ is Gorenstein under 
similar assumptions, this time the equality of ranks is not required.

The trivial source indecomposable $\calH_n$-modules are known as the 
\emph{Young modules} \cite{Dip-J89}. As shown in \cite{Dip-J91}, they are 
parametrized by partitions of $n$. Arbitrary indecomposable $\calH_n$-modules 
with a 1-dimensional source may be called \emph{signed Young modules} as in 
the case of representations of symmetric groups \cite{Don01}, \cite{Gia-LDW17}. 
However, an earlier text of Donkin \cite{Don} uses this term in a more 
restricted sense. There may be more such modules than partitions of $n$, and 
then the $\calH_n$-modules are not distinguished by their images in $\Sym$. 
Because of this we cannot generalize Theorem 6.1 to Hecke symmetries 
satisfying the 1-dimensional source condition.

In the semisimple case the algebras $A(R)$ and $A(R',R)$ also have rational 
Hilbert series. Moreover, \PHHai gives a formula for $H_{A(R'\!,R)}$ in 
terms of $H_{\bbS(V,R)}$ and $H_{\bbS(V'\!,R')}$ \cite{Hai02, Th. 3.1}. In 
section 5 these results are extended to the case of arbitrary $q$ under the 
assumption that both $R$ and $R'$ satisfy the $1$-dimensional source 
condition.

\section
1. Preliminaries

We fix an arbitrary field $\bbk$. Unless specified otherwise algebras and 
coalgebras will be considered over $\bbk$. Let $V$ be a finite dimensional 
vector space over $\bbk$. A \emph{Hecke symmetry} on $V$ is a linear operator 
$R:V\ot V\to V\ot V$ satisfying the braid equation
$$ 
(R\ot\Id)(\Id\ot\,R)(R\ot\Id)=(\Id\ot\,R)(R\ot\Id)(\Id\ot\,R) 
$$
and the quadratic Hecke relation
$$
(R-q\cdot\Id)(R+\Id)=0\quad\hbox{where $\,0\ne q\in\bbk$}.
$$

Denote by $\calH_n(q)$ the Hecke algebra of type $A_{n-1}$ with the same 
parameter $q$ as in the quadratic relation imposed on $R$. Since $R$ and $q$ 
will generally be fixed, we do not indicate $q$ in the notation $\calH_n$ 
when there is no danger of confusion. The algebra $\calH_n$ is generated by 
$n-1$ elements $T_1,\ldots,T_{n-1}$ subject to the defining relations
$$
\openup1\jot
\displaylines{
T_iT_j=T_jT_i\quad\hbox{whenever $|i-j|>1$},\cr
T_iT_{i+1}T_i=T_{i+1}T_iT_{i+1}\quad\hbox{for $i=1,\ldots,n-2$},\cr
(T_i-q)(T_i+1)=0\quad\hbox{for $i=1,\ldots,n-1$}.\cr
}
$$
Let $\frS_n$ be the symmetric group of permutations of the set $\{1,\ldots,n\}$. 
It is generated by basic transpositions $\tau_i=(i,i+1)$, $\,0<i<n$. Denote by 
$\ell(\si)$ the length of a permutation $\si\in\frS_n$ with respect to these 
generators. There is a standard basis $\{T_\si\mid\si\in\frS_n\}$ of $\calH_n$ 
characterized by the properties that $T_{\tau_i}=T_i$ for each $i$ and 
$T_{\pi\si}=T_{\pi}T_{\si}$ for $\pi,\si\in\frS_n$ whenever 
$\ell(\pi\si)=\ell(\pi)+\ell(\si)$. We adopt the convention that 
$\calH_0=\calH_1=\bbk$. 

The Hecke symmetry $R$ gives rise to a representation of $\calH_n$ in the 
$n$th tensor power of $V$ such that $T_i$ acts on $V^{\ot n}$ as the linear 
operator
$$
R_i^{(n)}=\Id^{\ot(i-1)}\!\ot\,R\ot\Id^{\ot(n-i-1)}\!. 
$$
In this way $V^{\ot n}$ becomes a left $\calH_n$-module.

Denote by $A(R)$ the \emph{$R$-matrix bialgebra}. It decomposes as a direct 
sum of subcoalgebras
$$
A(R)={\textstyle\bigoplus\limits_{n=0}^\infty}A_n(R)
$$
where $A_n(R)$ is the coalgebra dual to the finite dimensional algebra 
$\EndHn V^{\ot n}$.

Let $\Coend V$ be the coalgebra dual to $\EndkV$. For each $n\ge0$ we may 
identify $(\Coend V)^{\ot n}$ with the dual of the algebra 
$(\EndkV)^{\ot n}\cong\EndkV^{\ot n}$. Since $\EndHn V^{\ot n}$ is a 
subalgebra of $\EndkV^{\ot n}$, we have
$$
A_n(R)\cong(\Coend V)^{\ot n}/I_n
$$
where $I_n=(\EndHn V^{\ot n})^\perp=\{f\in(\Coend V)^{\ot n}\mid
\langle f,\,\EndHn V^{\ot n}\rangle=0\}$ is a coideal of 
$(\Coend V)^{\ot n}$. Clearly $I_n=0$ for $n=0,\mskip1mu 1$. For $n>1$ there 
is an equality
$$
\EndHn V^{\ot n}=\bigcap_{i=1}^{n-1}E_i^{(n)}
$$
where $E_i^{(n)}$ stands for the centralizer of $R_i^{(n)}$ in $\EndkV^{\ot n}$. 
Since
$$ 
E_i^{(n)}=(\EndkV)^{\ot(i-1)}\ot(\End_{\calH_2}\mskip-3mu 
V^{\ot2})\ot(\EndkV)^{\ot(n-i-1)}
$$
for each $i$, we get
$$
I_n=\sum_{i=1}^{n-1}{E_i^{(n)}}^{\,\perp}
=\sum_{i=1}^{n-1}\,(\Coend V)^{\ot(i-1)}\ot I_2\ot(\Coend V)^{\ot(n-i-1)}.
$$
This shows that $I=\bigoplus_{n=0}^\infty I_n$ is an ideal of the tensor algebra
$$
\bbT(\Coend V)={\textstyle\bigoplus\limits_{n=0}^\infty}\,(\Coend V)^{\ot n}
$$
generated by the homogeneous component $I_2$ of degree 2. By an earlier 
observation $I$ is also a coideal. Therefore $A(R)\cong\bbT(\Coend V)/I$ gets 
the structure of a factor bialgebra of $\bbT(\Coend V)$. This bialgebra 
coacts on $V$ universally with respect to the property that the induced 
coaction on $V^{\ot2}$ commutes with $R$. As observed in \cite{Lar-T91}, this 
property characterizes the bialgebra arising from the FRT construction.

For an associative algebra $\frA$ over some field we denote by $\Grot\frA$ the 
Grothendieck group of the category of finite dimensional left $\frA$-modules. 
To each finite dimensional left $\frA$-module $X$ there corresponds an element 
$[X]\in\Grot\frA$, and to each short exact sequence $0\to X'\to X\to X''\to0$ 
of finite dimensional left $\frA$-modules there corresponds a relation 
$[X]=[X']+[X'']$ in this group. The elements corresponding to finite dimensional 
left $\frA$-modules form a subsemigroup of $\Grot\frA$. Given $\xi\in\Grot\frA$, 
we write $\xi\ge0$ if $\xi$ lies in that subsemigroup, i.e., if $\xi=[X]$ for 
some finite dimensional left $\frA$-module $X$.

By the definition we have given above
$$
A_n(R)=\bigl(\EndHn V^{\ot n}\bigr)^*.
$$
Therefore right $A_n(R)$-comodules may be identified with left modules for the 
algebra $\EndHn V^{\ot n}$. The Grothendieck group $\Grot_n(R)$ of the category 
of finite dimensional right $A_n(R)$-comodules is identified with the group 
$\,\Grot(\EndHn V^{\ot n})$. The Grothendieck group of the category of finite 
dimensional right $A(R)$-comodules is the direct sum
$$ 
\Grot(R)={\textstyle\bigoplus\limits_{n=0}^\infty}\Grot_n(R).
$$
Moreover, $\Grot(R)$ is a graded ring with respect to the multiplication 
induced by tensor products of comodules.

The algebras $\bbS(V,R)$ and $\La(V,R)$ are defined as the factor algebras 
of the tensor algebra $\bbT(V)=\bigoplus_{n=0}^\infty\,V^{\ot n}$ by the 
ideals generated, respectively, by the subspaces
$$ 
\Im\,(R-q\cdot\Id)\sbs V^{\ot 2}\quad{\rm and}\quad
\Ker\,(R-q\cdot\Id)\sbs V^{\ot 2}.
$$
These ideals are stable under the coaction of $A(R)$ on $\bbT(V)$, and therefore 
$\bbS(V,R)$ and $\La(V,R)$ are right $A(R)$-comodule algebras in a natural way. 
The homogeneous components $\bbS_n(V,R)$ and $\La_n(V,R)$ of these algebras are 
right $A_n(R)$-comodules for each $n\ge0$. 

A \emph{composition} of $n$ is any finite sequences of positive integers 
$\la=(\la_1,\ldots,\la_k)$ with $|\la|=n$ where the \emph{weight} of $\la$ is 
defined as $|\la|=\sum\la_i$. The \emph{length} of $\la$ is the number 
$\ell(\la)=k$ of its parts $\la_i$. As is done customarily, we extend $\la$ by 
putting $\la_i=0$ for $i>\ell(\la)$. If $\la_i\ge\la_{i+1}$ for all $i$, 
then $\la$ is called a \emph{partition} of $n$. Denote by $\Par(n)$ the set of 
all partitions of $n$ and by $\Par$ the disjoint union of the sets 
$\Par(0),\,\Par(1),\,\Par(2),\ldots$ where $\Par(0)$ is regarded as a single 
element set consisting of the \emph{zero partition} $0$.

We denote by $\,\Sym\,$ the \emph{ring of symmetric functions} in a 
countable set of commuting indeterminates $x_1,x_2\mskip1mu,\ldots$ defined as 
in Macdonald \cite{Mac}. Thus, if $\,\Sym(r)\,$ stands for the subring of 
symmetric polynomials in the ring $\bbZ[x_1,\ldots,x_r]$, then $\Sym$ is the 
limit of the inverse system
$$
\cdots\to\Sym(3)\to\Sym(2)\to\Sym(1)\to\bbZ
$$
in the category of graded rings. If $u\in\Sym$ and 
$\al=(\al_1,\ldots,\al_r)\in K^r$ where $K$ is any commutative ring, then 
$u(\al)\in K$ is defined as the value at $\al$ of the polynomial obtained by 
projecting $u$ to $\Sym(r)$.

We conform to standard notation in regard to several families of symmetric 
functions \cite{Mac}. The elementary and complete symmetric functions will be 
$e_n$ and $h_n$ ($n=0,1,\ldots\,$) with $e_0=h_0=1$. For each partition 
$\la=(\la_1,\ldots,\la_k)$ one defines $e_\la=e_{\la_1}\!\cdots e_{\la_k}$ and 
$h_\la=h_{\la_1}\!\cdots h_{\la_k}$. The monomial functions $m_\la$ and the 
Schur functions $s_\la$ labelled by partitions $\la$ form two other well known 
$\bbZ$-bases of $\Sym$. On each homogeneous component $\Sym_n$ of the graded 
ring $\Sym$ there is a standard scalar product $\<\cdot\,,\cdot\>$ with 
respect to which $\{s_\la\mid\la\in\Par(n)\}$ is an orthonormal basis.

To each composition $\la$ of $n$ there corresponds a \emph{parabolic subalgebra} 
$\calH_\la=\calH_\la(q)$ of the Hecke algebra $\calH_n=\calH_n(q)$. It is 
generated by the set $\{T_i\mid i\in\calI_\la\}$ where
$$
\calI_\la=\{i\in\bbN\mid1\le i<|\la|\hbox{ and }i\ne\la_1+\ldots+\la_j
\hbox{ for each }j=1,\ldots,\ell(\la)\}
$$
and has a basis $\{T_\si\mid\si\in\frS_\la\}$ where $\frS_\la$ is the 
\emph{Young subgroup} of $\frS_n$ generated by $\{\tau_i\mid i\in\calI_\la\}$.

Each homomorphism $\chi:\calH_\la\to\bbk$ is completely determined by its 
values on the generators. It follows from the Hecke relations that 
$\chi(T_i)\in\{-1,q\}$ for each $i\in\calI_\la$. If both $i$ and $i+1$ lie in 
$\calI_\la$ then $\chi(T_i)=\chi(T_{i+1})$ by the braid relations. In other 
words, $\chi$ is constant on each of the $\ell(\la)$ contiguous (possibly empty) 
segments of lengths $\la_i-1$, $\,i=1,\ldots,\ell(\la)$, which comprise the 
set $\calI_\la$.

If $\nu$ is a composition of $n$ obtained from $\la$ by permuting its 
components in an arbitrary order, then the two subalgebras $\calH_\la$ and 
$\calH_\nu$ are conjugate by an inner automorphism of $\calH_n$. In this case 
the induction functors from $\calH_\la$ and from $\calH_\nu$ produce 
isomorphic $\calH_n$-modules. Given a homomorphism $\chi:\calH_\la\to\bbk$, it 
is possible to pass to an equivalent homomorphism $\chi':\calH_\nu\to\bbk$ 
such that all segments on which $\chi'$ takes value $q$ precede all segments 
on which $\chi'$ takes value $-1$, and any pair of segments on which $\chi'$ 
takes the same value follow in nonincreasing order of their lengths.

When forming the $\calH_n$-modules induced from 1-dimensional representations 
of parabolic subalgebras, it suffices to consider only the pairs 
$(\calH_\nu,\chi')$ satisfying the previous conditions. This leads to a 
parametrization of such modules by the set
$$ 
\Par^2(n)=\{(\la,\mu)\in\Par\times\Par\mid\,|\la|+|\mu|=n\}.
$$
Each pair $(\la,\mu)\in\Par^2(n)$ determines a composition 
$$
(\la_1,\ldots,\la_{\ell(\la)},\mu_1,\ldots,\mu_{\ell(\mu)}).
$$
We will denote by $\calH_\lamu$ the corresponding parabolic subalgebra of 
$\calH_n$, by $\frS_\lamu$ the corresponding Young subgroup of $\frS_n$, and 
by $\calI_\lamu$ the index set for the generators $T_i\in\calH_\lamu$ and 
$\tau_i\in\frS_\lamu$. Then
$$
\eqalign{
\calI_\lamu=\calI^0_\lamu\cup\calI^1_\lamu\quad{\rm where}\quad
\calI^0_\lamu&{}=\{i\in\calI_\lamu\mid i<|\la|\}=\calI_\la\,,\cr
\calI^1_\lamu&{}=\{i\in\calI_\lamu\mid i>|\la|\}
=\{|\la|+i\mid i\in\calI_\mu\}\,.
}
$$
Note that $\calH_\lamu\cong\calH_\la\ot\calH_\mu$. Define a homomorphism 
$\chi_\lamu:\calH_\lamu\to\bbk$ by the rule
$$
\chi_\lamu(T_i)=\cases{q & for $i\in\calI^0_\lamu\,$,\cr
\noalign{\smallskip}
-1 & for $i\in\calI^1_\lamu\,$,
}
$$
and denote by $\bbk_\lamu$ the corresponding 1-dimensional $\calH_\lamu$-module. 

If $\mu=0$, then $\calH_\lamu=\calH_\la$ and $\bbk_\lamu$ is the \emph{trivial} 
$\calH_\la$-module $\bbk\triv$ on which each generator $T_i\in\calH_\la$ 
operates as multiplication by $q$. If $\la=0$, then $\calH_\lamu=\calH_\mu$ 
and $\bbk_\lamu$ is the \emph{alternating} $\calH_\mu$-module $\bbk\alt$ on 
which each $T_i\in\calH_\mu$ operates as multiplication by $-1$.

From the preceding discussion it follows that an indecomposable left 
$\calH_n$-module has a \emph{$1$-dimensional} (respectively, \emph{trivial}) 
\emph{source} if and only if it is isomorphic to a direct summand of the 
induced module $\calH_n\ot_{\calH_{\la,\mu}}\bbk_\lamu$ (respectively, 
$\calH_n\ot_{\calH_\la}\bbk\triv$) for some $(\la,\mu)\in\Par^2(n)$ 
(respectively, $\la\in\Par(n)$).

The Specht $\calH_n$-modules $S^\la$ labelled by partitions $\la\in\Par(n)$ 
were constructed by Dipper and James \cite{Dip-J86}. We will use this notation 
for left modules and sometimes also for right modules as in \cite{Dip-J86}. 
The dimension of $S^\la$ depends neither on the field $\bbk$ nor on the 
parameter $q$. In particular, it is the same as the dimension of the 
respective Specht module for the symmetric group $\frS_n$.

The Hecke algebras $\calH_0,\calH_1,\calH_2,\ldots$ are all semisimple if and 
only if $[n]_q\ne0$ for all $n>0$ where
$$
[n]_q=1+q+\cdots+q^{n-1}.
$$
We now recall briefly what happens in the \emph{semisimple case}. The simple 
$\calH_n$-modules are precisely the Specht modules. Moreover, each $S^\la$ is 
absolutely irreducible, so that $\End_{\calH_n}\!S^\la\cong\bbk$. By 
semisimplicity of $\calH_n$ an arbitrary $\calH_n$-module $N$ is a direct sum 
of its isotypic components, and the isotypic component of type $S^\la$ can be 
expressed as $\Hom_{\calH_n}(S^\la\!,N)\ot S^\la$. Taking $N=V^{\ot n}\!$, we get
$$
V^{\ot n}\cong{\textstyle\bigoplus\limits_{\la\in\Par(n)}}V_R^\la\ot S^\la\quad 
\hbox{where }V_R^\la=\Hom_{\calH_n}(S^\la,V^{\ot n}).
$$
It follows from this decomposition that
$$
\EndHn V^{\ot n}\cong\prod_{\la\in\Par(n)}\EndkV_R^\la.
$$
Thus the algebra $\EndHn V^{\ot n}$ is semisimple with 
$\{V_R^\la\mid\la\in\Par(n),\ V_R^\la\ne0\}$ being a full set of pairwise 
nonisomorphic simple left modules. Note that $V_R^\la\ne0$ if and only if 
$S^\la$ embeds in $V^{\ot n}$ as an $\calH_n$-submodule. We see that all right 
$A(R)$-comodules are semisimple, which means that $A(R)$ is a 
\emph{cosemisimple} bialgebra. Moreover,
$$
\{V_R^\la\mid\la\in\Par,\ V_R^\la\ne0\}
$$
is a full set of pairwise nonisomorphic simple right $A(R)$-comodules. Their 
isomorphism classes $[V^\la]$ form a $\bbZ$-basis of the Grothendieck group 
$\Grot(R)$. Note also that $\dim S^\la$ is the multiplicity of $V_R^\la$ as an 
$A(R)$-comodule summand of $V^{\ot n}$, while $\dim V_R^\la$ is the 
multiplicity of $S^\la$ as an $\calH_n$-module summand of $V^{\ot n}$.

Let $\calH_{m,n}$ be the subalgebra of $\calH_{m+n}$ generated by 
$\{T_i\mid0<i<m+n,\ i\ne\nobreak m\}$. In other words, 
$\calH_{m,n}=\calH_{(m),(n)}$ for partitions $(m),(n)$ of length 1.
Denote by $\ind_{m,n}^{m+n}$ the induction functor from $\calH_{m,n}$ to 
$\calH_{m+n}$. There is a graded ring structure on the direct sum of 
Grothendieck groups $\bigoplus_{k=0}^\infty\Grot\calH_k$ defined by the rule
$$
[X]\cdot[\mskip1mu Y]=[\ind_{m,n}^{m+n}X\ot Y]
$$
whenever $X$ is an $\calH_m$-module and $Y$ an $\calH_n$-module, both of 
finite dimension. Here $X\ot Y$ is viewed as an $\calH_{m,n}$-module by means 
of the canonical isomorphism $\calH_{m,n}\cong\calH_m\ot\calH_n$. In this ring 
$$
[S^\mu]\cdot[S^\nu]=\sum_{\la\in\Par}c_{\mu\nu}^\la\,[S^\la]\quad
\hbox{for $\mu,\nu\in\Par$}
$$
with the Littlewood-Richardson coefficients $c_{\mu\nu}^\la$ which also 
occur as structure constants for the multiplication in the ring of symmetric 
functions:
$$
s_\mu s_\nu=\sum_{\la\in\Par}c_{\mu\nu}^\la\,s_\la.
$$
This means that there is an isomorphism of graded rings
$$
\Sym\cong{\textstyle\bigoplus\limits_{k=0}^\infty}\,\Grot\calH_k
$$
under which $s_\la\in\Sym$ corresponds to $[S^\la]\in\Grot\calH_{|\la|}$. 
Conceptual explanation of this isomorphism is provided by Zelevinsky's 
approach \cite{Zel}. With some additional structure the direct sum of the 
groups $\Grot\calH_k$ satisfies the axioms of a connected positive selfadjoint 
Hopf algebra over $\bbZ$, and it contains only one irreducible primitive 
element. It was proved in \cite{Zel} that such a Hopf algebra is unique up to 
isomorphism and is isomorphic to the ring of symmetric functions.

\proclaim
Lemma 1.1.
For each $(\la,\mu)\in\Par^2(n)$ the symmetric function $h_\la e_\mu$ maps to 
the class of the induced module $\calH_n\ot_{\calH_{\la,\mu}}\bbk_\lamu$ in the 
Grothendieck group\/ $\Grot\calH_n$. As a consequence{\rm,} the following 
relations hold in this group\/{\rm:}
$$
\displaylines{
\sum_{i=0}^n\,(-1)^i\,[\calH_n\ot_{\calH_{i,n-i}}\bbk\iin]=0\quad
{\rm if}\ n>0,\cr
[\calH_n\ot_{\calH_\mu}\bbk\triv]=\sum_{\la\in\Par(n)}K_{\la\mu}\,[S^\la]\quad
{\rm for}\ \mu\in\Par(n)\cr
}
$$
where $\bbk\iin$ is the $1$-dimensional $\calH\iin$-module associated with 
the pair of partitions $\bigl((i),(n-i)\bigr)\in\Par^2(n)$ and $K_{\la\mu}$ are 
the Kostka numbers.
\endproclaim

\Proof.
Since for each $p>0$ the symmetric functions $h_p=s_{(p)}$ and $e_p=s_{(1^p)}$ 
map to the classes of the Specht modules $S^{(p)}=\bbk\triv$ and 
$S^{(1^p)}=\bbk\alt$ in the group $\Grot\calH_p$, it follows from the 
definition of the multiplication in the direct sum of the groups $\Grot\calH_k$ 
that $h_\la$ and $e_\la$, for any $\la\in\Par$, map to the classes of induced 
modules $\calH_{|\la|}\ot_{\calH_\la}\bbk\triv$ and 
$\calH_{|\la|}\ot_{\calH_\la}\bbk\alt$ in the group $\Grot\calH_{|\la|}$.
Hence
$$
h_\la e_\mu\mapsto[\calH_{|\la|}\ot_{\calH_\la}\bbk\triv]\cdot
[\calH_{|\mu|}\ot_{\calH_\mu}\bbk\alt]=[\calH_n\ot_{\calH_{\la,\mu}}\bbk_\lamu]
$$
for $(\la,\mu)\in\Par^2(n)$. The required equalities in the group 
$\Grot\calH_n$ are now immediate consequences of the well known equalities 
$$
\sum_{i=0}^n(-1)^ih_ie_{n-i}=0,\qquad\quad
h_\mu=\sum_{\la\in\Par(n)}K_{\la\mu}s_\la
$$
in the group $\Sym_n$ (see \cite{Mac, Ch. I}).
\endproof

In the semisimple case tensor products of simple $A(R)$-comodules are computed 
easily. If $\mu\in\Par(m)$ and $\nu\in\Par(n)$, then
$$
\openup1\jot
\eqalign{
V_R^\mu\ot V_R^\nu&{}=\Hom_{\calH_m}(S^\mu,V^{\ot m})\ot
\Hom_{\calH_n}(S^\nu,V^{\ot n})\cr 
&{}\cong\Hom_{\calH_{m,n}}(S^\mu\ot S^\nu,\,V^{\ot m}\ot V^{\ot n})\cr 
&{}\cong\Hom_{\calH_{m+n}}(\ind_{m,n}^{m+n}S^\mu\ot S^\nu,\,V^{\ot(m+n)})\,.
}
$$
It follows that for each $\la\in\Par(m+n)$ the multiplicity of $V_R^\la$ in 
$V_R^\mu\ot V_R^\nu$ equals the multiplicity of the simple $\calH_{m+n}$-module 
$S^\la$ in $\,\ind_{m,n}^{m+n}S^\mu\ot S^\nu$. Thus
$$
[V_R^\mu]\cdot[V_R^\nu]=\sum_{\la\in\Par}\,c_{\mu\nu}^\la\,[V_R^\la]\quad
\hbox{in $\Grot(R)$},
$$
and so there is a surjective ring homomorphism $\ph:\Sym\to\Grot(R)$ given by 
the assignments $s_\la\mapsto[V_R^\la]$, $\,\la\in\Par$. Let
$$
\Ga(r_0,r_1)=\{\la\in\Par\mid\la_j\le r_1\hbox{ for all }j>r_0\}
$$
be the set of $(r_0,r_1)$-hook partitions. It was proved by \PHHai 
\cite{Hai99} that $V_R^\la\ne0$ if and only if $\la\in\Ga(r_0,r_1)$ where 
$(r_0,r_1)$ is the birank of $R$. Hence the classes $[V_R^\la]$ with 
$\la\in\Ga(r_0,r_1)$ form a $\bbZ$-basis of $\Grot(R)$ and $\Ker\ph$ coincides 
with the $\bbZ$-linear span of $\{s_\la\mid\la\notin\Ga(r_0,r_1)\}$.

In the present paper it will be shown that some features of this situation 
extend to the \emph{nonsemisimple case} provided that $R$ satisfies the 
1-dimensional source condition.

\section
2. The decomposition map

The decomposition map is a standard tool in the modular representation theory 
of finite groups. More generally, such a map can be defined in the following 
situation (see, e.g., \cite{Geck-P, 7.4.3}). Suppose that $O$ is a discrete 
valuation ring with residue field $\bbk$ and the field of fractions $Q$. Let 
$A$ be an associative unital algebra over $O$ whose underlying $O$-module is 
free of finite rank. Then $\Ak=A\ot_O\bbk$ and $A_Q=A\ot_OQ$ are finite 
dimensional algebras over the respective fields. By an \emph{$A$-lattice} we 
mean any finitely generated $O$-free $A$-module (a finitely generated 
$O$-module is free if and only if it is torsionfree). The decomposition map
$$
d:\Grot A_Q\mapr{}\Grot\Ak
$$
is a homomorphism of groups characterized by the property that 
$d([L_Q])=[\Lk\mskip0.5mu]$ 
for each left $A$-lattice $L$ where $L_Q=L\ot_OQ$ and $\Lk=L\ot_O\bbk$. This 
map is well-defined since each $A_Q$-module of finite dimension over $Q$ is 
isomorphic to $L_Q$ for some $A$-lattice $L$, and the image of $\Lk$ in 
$\Grot\Ak$ does not depend on the choice of $L$.

Let $z$ be an invertible element of $O$ and $\calH_n(z)$ the Hecke algebra of 
type $A_{n-1}$ with parameter $z$ over the ring $O$. If $q$ is the image of 
$z$ in $\bbk$, then $\calH_n(z)_\bbk\cong\calH_n(q)$. We will assume that
$$ 
1+z+\ldots+z^{i-1}\ne0\quad\hbox{for all $i>0$}.
$$
This can be achieved, e.g., by taking $z$ to be an indeterminate and $O$ the 
localization of the polynomial ring $\bbk[z]$ at its maximal ideal generated 
by $z-q$. Then $\calH_n(z)_Q$ is a semisimple Hecke algebra of type $A_{n-1}$ 
over the field $Q$. By completing $O$ we may also assume that $O$ is a 
\emph{complete discrete valuation ring}.

We will denote by $T_1,\ldots,T_{n-1}$ the canonical generators of $\calH_n(z)$ 
and also their canonical images in either $\calH_n(q)$ or $\calH_n(z)_Q$. 
For $(\la,\mu)\in\Par^2(n)$ we have a parabolic subalgebra $\calH_\lamu(z)$ 
generated by $\{T_i\mid i\in\calI_\lamu\}$ (see section 1). Define a 
homomorphism of $O$-algebras $\chi_\lamu:\calH_\lamu(z)\to O$ by the rule
$$
\chi_\lamu(T_i)=\cases{z & for $i\in\calI^0_\lamu\,$,\cr
\noalign{\smallskip}
-1 & for $i\in\calI^1_\lamu\,$,
}
$$
and let $O_\lamu$ be $O$ with the $\calH_\lamu(z)$-module structure arising 
from $\chi_\lamu$. Let $Q_\lamu$ be the 1-dimensional 
$\calH_\lamu(z)_Q$-module defined similarly. We will write $\calH_\la(z)$, 
$O\triv$, and $Q\triv$ instead of $\calH_\lamu(z)$, $O_\lamu$, and 
$Q_\lamu$ when $\la\in\Par(n)$ and $\mu=0$.

\medbreak
{\bf Notation.}
Denote by $\Rep$ (respectively, by $\Triv$) the class of all 
$\calH_n(z)$-lattices isomorphic to direct summands of finite direct sums of 
the induced modules
$$
M^\lamu=\calH_n(z)\ot_{\calH_{\la,\mu}(z)}O_\lamu\quad 
{\rm(respectively,}\quad M^\la=\calH_n(z)\ot_{\calH_\la(z)}O\triv{\rm)}
$$
for various $(\la,\mu)\in\Par^2(n)$ (respectively, $\la\in\Par(n)$).
\medbreak

Note that $M^\la=M^\lazero$. Since $\calH_n(z)$ is a free 
$\calH_\lamu(z)$-module with respect to the action by right 
multiplications, $M^\lamu$ is an $\calH_n(z)$-lattice. The functor 
$?\ot_O\bbk$ takes $M^\lamu$ to
$$
\Mk^\lamu\cong\calH_n(q)\ot_{\calH_{\la,\mu}(q)}\bbk_\lamu\,.
$$
It follows that all indecomposable direct summands of the $\calH_n(q)$-module 
$\Mk$ have a 1-dimensional (respectively, trivial) source whenever 
$M\in\Rep$ (respectively, $M\in\Triv$).

\proclaim
Lemma 2.1.
Let $M,N\in\Rep$. If $q=-1$ assume that $M,N\in\Triv$. Then
$$
\openup1\jot
\eqalign{
\Hom_{\calH_n(q)}(\Nk,\Mk)&{}\cong\Hom_{\calH_n(z)}(N,M)\ot_O\bbk\,,\cr
\Hom_{\calH_n(z)_Q}(N_Q,M_Q)&{}\cong\Hom_{\calH_n(z)}(N,M)\ot_OQ.
}
$$
Hence\quad$\dimk\Hom_{\calH_n(q)}(\Nk,\Mk)=\dim_Q\Hom_{\calH_n(z)_Q}(N_Q,M_Q)$.
\endproclaim

\Proof.
The canonical maps 
$k_{XY}:\Hom_{\calH_n(z)}(X,Y)\ot_O\bbk\to\Hom_{\calH_n(q)}(X_\bbk,Y_\bbk)$ 
defined for arbitrary $\calH_n(z)$-modules $X$ and $Y$ give a natural 
transformation of two functors additive in each argument. If $X\cong X'\oplus X''$ 
(respectively, $Y\cong Y'\oplus Y''$), then bijectivity of $k_{XY}$ is 
equivalent to bijectivity of $k_{X'Y}$ and $k_{X''Y}$ (respectively, 
$k_{XY'}$ and $k_{XY''}$).
By the conditions on $M$ and $N$ in the statement of Lemma 2.1 it suffices 
therefore to prove that $k_{NM}$ is bijective when $q\ne-1$ and 
$N=M^\lamu$, $M=M^\rhth$ for some pairs $(\la,\mu),\,(\rho,\th)\in\Par^2(n)$ 
or when $q=-1$ and $N=M^\la$, $M=M^\rho$ with $\la,\rho\in\Par(n)$.

Consider the case $q\ne-1$. Denote by $\calD$ the set of distinguished 
representatives of the $\frS_\lamu\,$-$\,\frS_\rhth$ double cosets in $\frS_n$. 
Then $M=\bigoplus_{\pi\in\calD}\,M(\pi)$ where $M(\pi)$ is the 
$\calH_\lamu(z)$-submodule of $M$ generated by the element $T_\pi\ot1\in M$. 
This is the Mackey decomposition of the induced module $M$ with respect to the 
parabolic subalgebra $\calH_\lamu(z)$ (see \cite{Dip-J86, 2.7} and 
\cite{Geck-P, 9.1.8}). For each $\pi\in\calD$ let $\nu(\pi)$ be the 
composition of $n$ such that
$$
\eqalign{
\calI_{\nu(\pi)}&{}=\{i\in\calI_\lamu\mid\pi^{-1}\tau_i\pi\in\frS_\rhth\}\cr
&{}=\{i\in\calI_\lamu\mid\pi^{-1}\tau_i\pi=\tau_j\hbox{ for some }
j\in\calI_\rhth\}.
}
$$
If $i\in\calI_{\nu(\pi)}$, then $\ell(\tau_i\pi)>\ell(\pi)$ since $\pi$ is 
the shortest element in the coset $\frS_\lamu\,\pi$, and the equality 
$\tau_i=\pi\mskip1mu\tau_j\pi^{-1}$ means that $i=\pi(j)$ and $i+1=\pi(j+1)$. 
Therefore $T_iT_\pi=T_{\tau_i\pi}=T_\pi T_j$ with $j=\pi^{-1}(i)\in\calI_\rhth$, 
i.e., $T_\pi^{-1}T_iT_\pi=T_{\pi^{-1}(i)}\in\calH_\rhth(z)$. 
By the Mackey formula
$$
M(\pi)\cong\calH_\lamu(z)\ot_{\calH_{\nu(\pi)}(z)}O(\chi'_\pi)
$$
where $O(\chi'_\pi)$ is $O$ regarded as an $\calH_{\nu(\pi)}(z)$-module by 
means of the ring homomorphism $\chi'_\pi:\calH_{\nu(\pi)}(z)\to O$ such that 
$$
\chi'_\pi(h)=\chi_\rhth(T_\pi^{-1}hT_\pi)\quad
\hbox{for $h\in\calH_{\nu(\pi)}(z)$}.
$$
Note that
$$
\chi'_\pi(T_i)=\chi_\rhth(T_{\pi^{-1}(i)})
=\cases{z & for $i\in\calI_{\nu(\pi)}$ with 
$\pi^{-1}(i)\in\calI_\rhth^0\,$,\cr
\noalign{\smallskip}
-1 & for $i\in\calI_{\nu(\pi)}$ with $\pi^{-1}(i)\in\calI_\rhth^1\,$.
}
$$
By the Frobenius reciprocity (see \cite{Dip-J86, 2.5, 2.6} and 
\cite{Geck-P, 9.1.7})
$$
\Hom_{\calH_n(z)}(N,M)\cong\Hom_{\calH_{\la,\mu}(z)}(O_\lamu,M) 
\cong\bigoplus_{\pi\in\calD}\,\Hom_{\calH_{\nu(\pi)}(z)}\bigl(O_\lamu,O(\chi'_\pi)\bigr).
$$
The $O$-module $\Hom_{\calH_{\nu(\pi)}(z)}\bigl(O_\lamu,O(\chi'_\pi)\bigr)$ 
is nonzero precisely when $\chi_\lamu$ agrees with $\chi'_\pi$ on 
$\calH_{\nu(\pi)}(z)$, in which case this module is isomorphic to $O$. It 
follows that $\Hom_{\calH_n(z)}(M,N)$ is a free $O$-module with a basis 
indexed by the set of those $\pi\in\calD$ for which 
$\,\calI_{\nu(\pi)}\cap\calI_\lamu^0=\calI_{\nu(\pi)}\cap\pi(\calI_\rhth^0)$.

The respective homomorphisms $N\to M$ can be described explicitly as in 
\cite{Dip-J86, 3.4}. The vector space $\Hom_{\calH_n(q)}(\Nk,\Mk)$ has a 
similar description, and we see that the functor $?\ot_O\bbk$ takes the basic 
homomorphisms $N\to M$ to the basic homomorphisms $\Nk\to \Mk$. If 
$\mu=0$ and $\th=0$, then the basic homomorphisms are parametrized by the 
whole set $\calD$, and the previous arguments go through for $q=-1$ as well. 
This proves bijectivity of $k_{NM}$.

The second isomorphism in the statement of Lemma 2.1 can be explained in 
exactly the same way, but in fact it is almost obvious and holds more generally 
for arbitrary $\calH_n(z)$-modules $M$ and $N$ with the only restriction that 
$N$ should be finitely generated. In the last assertion of Lemma 2.1 both 
dimensions are equal to the rank of the free $O$-module $\Hom_{\calH_n(z)}(N,M)$.
\endproof

For the rest of this section with the exception of Lemma 2.7 we fix $M\in\Rep$, 
and moreover we will assume that $M\in\Triv$ when $q=-1$. The ring $A=\EndHM$ 
is an algebra over $O$ whose underlying $O$-module is free of finite rank. By 
Lemma 2.1 $\Ak\cong\EndHMk$ and $A_Q\cong\EndHMQ$. Thus we have the 
decomposition map
$$
d:\Grot(\EndHMQ)\mapr{}\Grot(\EndHMk).
$$

\proclaim
Lemma 2.2.
If $q\ne-1$ and $(\la,\mu)\in\Par^2(n)$ then 
$$
d(\,[\,Q_\lamu\ot_{\calH_{\la,\mu}(z)_Q}\!M_Q\,]\,)
=[\,\bbk_\lamu\ot_{\calH_{\la,\mu}(q)}\Mk\,]\,.
$$
In particular{\rm,} $\,d(\,[\,Q\triv\ot_{\calH_\la(z)_Q}\!M_Q\,]\,)
=[\,\bbk\triv\ot_{\calH_\la(q)}\Mk\,]\,$ for $\la\in\Par(n),$ and this 
equality holds even when $q=-1$.
\endproclaim

\Proof.
Put $L=O_\lamu\ot_{\calH_{\la,\mu}(z)}M$. This is an 
$\EndHM$-module such that 
$$
\Lk\cong\bbk_\lamu\ot_{\calH_{\la,\mu}(q)}\Mk\quad{\rm and}\quad 
L_Q\cong Q_\lamu\ot_{\calH_{\la,\mu}(z)_Q}\!M_Q.
$$
We will check that $L$ is $O$-free of finite rank, i.e. $L$ is a lattice. It 
will follow then that $d([L_Q])=[\Lk]$ by the definition of $d$, and we 
will get the required equality.

Since the verifications can be done on direct summands, it suffices to 
consider the case when $M=M^\rhth$ for some $(\rho,\th)\in\Par^2(n)$. Using 
the Mackey decomposition $M=\bigoplus_{\pi\in\calD}M(\pi)$ with respect to 
$\calH_\lamu(z)$ as in the proof of Lemma 2.1, we get
$L=\bigoplus_{\pi\in\calD}L(\pi)$ with
$$
L(\pi)=O_\lamu\ot_{\calH_{\la,\mu}(z)}M(\pi)\cong 
O_\lamu\ot_{\calH_{\nu(\pi)}(z)}O(\chi'_\pi)\cong O/I_\pi
$$
where $I_\pi$ is the ideal of $O$ generated by 
$\{\chi_\lamu(h)-\chi'_\pi(h)\mid h\in\calH_{\nu(\pi)}(z)\}$. Since the 
$O$-algebra $\calH_{\nu(\pi)}(z)$ is generated by 
$\{T_i\mid i\in\calI_{\nu(\pi)}\}$, the ideal $I_\pi$ is generated by the 
elements
$$
\chi_\lamu(T_i)-\chi'_\pi(T_i)\quad\hbox{with $\,i\in\calI_{\nu(\pi)}$}\,.
$$
If $\chi_\lamu$ agrees with $\chi'_\pi$ on $\calH_{\nu(\pi)}(z)$, then 
$I_\pi=0$. Otherwise $\chi_\lamu(T_i)\ne\chi'_\pi(T_i)$ for at least one 
$i\in\calI_{\nu(\pi)}$. Since $-1$ and $z$ are the only two possible values of 
the homomorphisms $\chi_\lamu$ and $\chi'_\pi$ on the generators $T_i$, we get 
$I_\pi=(z+1)O$ in the latter case. The image of $z+1$ in the residue field 
$\bbk$ of the local ring $O$ equals $q+1$. If $q\ne-1$, then $z+1$ is 
invertible in $O$, and so $I_\pi=O$. We see that $O/I_\pi$ equals either $O$ 
or $0$ for each $\pi\in\calD$ when $q\ne-1$, i.e. each $O$-module in the direct 
sum decomposition of $L$ is free of rank 1 or 0.

Suppose now that $M\in\Triv$. In this case we may assume that $M=M^\rho$ for 
some $\rho\in\Par(n)$, i.e. we can take $\th=0$. Then $\chi'_\pi$ is the 
trivial representation. If $\mu=0$, then $\chi_\lamu$ is the trivial 
representation as well, whence $I_\pi=0$ and $O/I_\pi\cong O$ for all $\pi$, 
even when $q=-1$.
\endproof

\medbreak
{\bf Notation.}
For any $\calH_n(q)$-module $X$ and $(\la,\mu)\in\Par^2(n)$ put
$$
\Si_\lamu(X)=\sum_{i\in\calI^0_\lamu}(T_i-q)X
+\sum_{i\in\calI^1_\lamu}\{x\in X\mid T_ix=qx\}.
$$
This is an $\EndHX$-submodule of $X$. In particular, we write
$$
\Si\iin(X)=\sum_{0<j<i}(T_j-q)X
+\sum_{i<j<n}\{x\in X\mid T_jx=qx\}
$$
when $\la=(i)$, $\mu=(n-i)$ for some $i=0,\ldots,n$, and we write
$$
\Si_\la(X)=\Si_{\la,0}(X)=\sum_{i\in\calI_\la}(T_i-q)X
$$
when $\la\in\Par(n)$, $\mu=0$.

\medbreak
The case $q=-1$ incurs technical complications, and we have to look deeper 
into the structure of induced modules. The conclusion of Lemma 2.2 is 
reformulated below in a form suitable for any $q$:

\proclaim
Lemma 2.3.
We have\quad 
$d(\,[\,Q_\lamu\ot_{\calH_{\la,\mu}(z)_Q}\!M_Q\,]\,)=[\,\Mk/\Si_\lamu(\Mk)\,]$.
\endproclaim

\Proof.
Put $L=O_\lamu\ot_{\calH_{\la,\mu}(z)}M$ and $A=\EndHM$. The assignment 
$m\mapsto1\ot m$ defines an epimorphism of $A$-modules $\psi:M\to L$. The 
induced map $\psik=\psi\ot_O\bbk$ is an epimorphism of $\Ak$-modules 
$\,\Mk\to \Lk\cong\bbk_\lamu\ot_{\calH_{\la,\mu}(q)}\Mk\,$ such that
$$
\eqalign{
\Ker\psik=\sum_{h\in\calH_\lamu(q)}\bigl(h-\chilamu(h)\bigr)\Mk
&{}=\sum_{i\in\calI_\lamu}\bigl(T_i-\chilamu(T_i)\bigr)\Mk\cr
&{}=\sum_{i\in\calI^0_\lamu}(T_i-q)\Mk+\sum_{i\in\calI^1_\lamu}(T_i+1)\Mk\,.
}
$$
Suppose that $q\ne-1$. It follows then from the relation $(T_i-q)(T_i+1)=0$ 
that the equality $T_im=qm$ holds for an element $m\in \Mk$ if and only 
if $m\in(T_i+1)\Mk$. Hence $\,\Ker\psik=\Si_\lamu(\Mk)$, and
$$ 
\Lk\cong \Mk/\Ker\psik=\Mk/\Si_\lamu(\Mk).
$$
The equality $d([L_Q])=[\Lk]$ of Lemma 2.2 can be rewritten as in the 
statement of Lemma 2.3.

Suppose further that $q=-1$. In this case $L$ may fail to be $O$-free. The set 
$L^0$ of all elements of $L$ which have nonzero annihilators in $O$ is an 
$A$-submodule of $L$. The factor module $\Lov=L/L^0$ is $O$-torsionfree and 
$\Lov_Q\cong L_Q$. Hence $\Lov$ is a lattice, which yields 
$\,d([L_Q])=[\mkern1mu\Lov_\bbk\mkern1mu]$.

We have $\Lov\cong M/M^0$ where $M^0=\psi^{-1}(L^0)\sbs M$. Consider 
$\Mk^0=M^0\ot_O\bbk$. The canonical map $\Mk^0\to\Mk$ is injective since 
$\Lov$ is $O$-free. Thus $\Mk^0$ is identified with an $\Ak$-submodule of 
$\Mk$. By the right exactness of the functor $?\ot_O\bbk$ we obtain 
$\Lov_\bbk\cong \Mk/\Mk^0$, and it remains to show that $\Mk^0=\Si_\lamu(\Mk)$.

Note that $\Mk^0$ and $\Si_\lamu(\Mk)$ are evaluations at $M$ of two 
additive functors defined on the category of $\calH_n(z)$-modules. Verification 
of the equality $\Mk^0=\Si_\lamu(\Mk)$ can be done therefore on direct 
summands of $M$. Since $M\in\Triv$, it suffices to consider the case when 
$M=M^\rho$ for some $\rho\in\Par(n)$. Let us assume this and take $\th=0$ in 
the proof of Lemma 2.2. We have seen there that $L=\bigoplus_{\pi\in\calD}L(\pi)$ 
with $L(\pi)\cong O/I_\pi$.

Since $\chi'_\pi$ is now the trivial representation for each $\pi\in\calD$, 
we have $I_\pi\ne0$ if and only if $\chi_\lamu(T_i)=-1$ for at least one 
$i\in\calI_{\nu(\pi)}$. This condition on $\pi$ means precisely that 
$\,\calI_{\nu(\pi)}\cap\calI^1_\lamu\ne\varnothing$. Hence
$$
L^0=\bigoplus_{\pi\in\calD,\ \calI_{\nu(\pi)}\cap\calI^1_\lamu\ne\varnothing}L(\pi).
$$
For each $\si\in\frS_n$ denote by $v_\si$ the element $T_\si\ot1\in 
M=\calH_n(z)\ot_{\calH_\rho(z)}O\triv$ and by $v'_\si$ a similar element of 
$\Mk$. Then $v'_\si$ is the image of $v_\si$ under the canonical map 
$M\to\Mk$. Since $L(\pi)$ is the cyclic $O$-submodule of $L$ generated by 
$\psi(v_\pi)$, we get
$$
\eqalign{
M^0&{}=\sum_{\pi\in\calD,\ \calI_{\nu(\pi)}\cap\calI^1_\lamu\ne\varnothing}Ov_\pi
+\Ker\psi\,,\qquad\hbox{and therefore}\cr
\Mk^0&{}=\sum_{\pi\in\calD,\ \calI_{\nu(\pi)}\cap\calI^1_\lamu\ne\varnothing}
\bbk v'_\pi+\Ker\psik\,.
}
$$
Put $U_i=\{m\in\Mk\mid T_im=-m\}$. Since $(T_i+1)^2=0$ in the algebra 
$\calH_n(q)$, we have $(T_i+1)\Mk\sbs U_i$. Hence
$$
\Ker\psik=\sum_{i\in\calI_\lamu}(T_i+1)\Mk
\sbs\sum_{i\in\calI_\lamu^0}(T_i+1)\Mk+\sum_{i\in\calI_\lamu^1}U_i
=\Si_\lamu(\Mk).
$$
If $\pi\in\calD$, then $T_iv'_\pi=-v'_\pi$, i.e. $v'_\pi\in U_i$, for each 
$i\in\calI_{\nu(\pi)}$. This shows that $v'_\pi$ lies in $\Si_\lamu(\Mk)$ 
whenever $\calI_{\nu(\pi)}\cap\calI^1_\lamu\ne\varnothing$. Hence 
$\Mk^0\sbs\Si_\lamu(\Mk)$.

Conversely, we claim that $U_i\sbs\Mk^0$ for each $i\in\calI^1_\lamu$, 
which entails the opposite inclusion $\Si_\lamu(\Mk)\sbs\Mk^0$. Fix 
such an $i$ and consider the Mackey decomposition 
$\Mk=\bigoplus_{\pi\in\calD}M(\pi)_\bbk$ where $M(\pi)_\bbk$ is the 
$\calH_\lamu(q)$-submodule of $\Mk$ generated by $v'_\pi$. The standard 
basis of $M(\pi)_\bbk\cong\calH_\lamu(q)\ot_{\calH_{\nu(\pi)}(q)}\bbk\triv$ is 
formed by the elements $v'_{\si\pi}=T_\si v'_\pi$ with $\si$ in the set 
$\calD^\lamu_\pi$ of distinguished representatives of the cosets 
$\frS_\lamu/\frS_{\nu(\pi)}$.

For each $\si\in\calD^\lamu_\pi$ either $\tau_i\si\in\calD^\lamu_\pi$ or, by 
Deodhar's Lemma, $\tau_i\si=\si\tau_j$ for some $j\in\calI_{\nu(\pi)}$. In 
the first case $v'_{\si\pi}$ and $v'_{\tau_i\si\pi}$ span a 2-dimensional 
$T_i$-invariant subspace whose intersection with $U_i$ is spanned by a single 
element 
$$
v'_{\tau_i\si\pi}+v'_{\si\pi}=\pm(T_i+1)v'_{\si\pi}\in\Ker\psik\sbs\Mk^0.
$$
In the second case $T_iT_\si=T_{\tau_i\si}=T_\si T_j$ and $T_jv'_\pi=-v'_\pi$, 
whence 
$$
T_iv'_{\si\pi}=T_iT_\si v'_\pi=T_\si T_jv'_\pi=-T_\si v'_\pi=-v'_{\si\pi}\,.
$$
Thus $v'_{\si\pi}$ spans a 1-dimensional $T_i$-invariant subspace. 
Note that $j\in\calI^1_\lamu$ since the transpositions $\tau_i$ and $\tau_j$ 
are conjugate in the group $\frS_\lamu$. Thus 
$\calI_{\nu(\pi)}\cap\calI^1_\lamu\ne\varnothing$, which entails 
$v'_\pi\in\Mk^0$. But $T_\si v'_\pi\equiv(-1)^{\ell(\si)}v'_\pi$ modulo 
$\Ker\psik$ since $\si\in\frS_\lamu$ and $(T_l+1)\Mk\sbs\Ker\psik$ 
for all $l\in\calI_\lamu$. Therefore $v'_{\si\pi}\in\Mk^0$ too.

The whole $\Mk$ is thus a direct sum of $T_i$-invariant subspaces spanned 
by at most 2 basis elements. It follows that $U_i$ is a direct sum of its 
intersections with those subspaces of $\Mk$. We have checked that $\Mk^0$ 
contains each of the summands in this decomposition of $U_i$. Hence 
$U_i\sbs\Mk^0$, as claimed.
\endproof

\proclaim
Lemma 2.4.
The equality $\ \sum\limits_{i=0}^n(-1)^i\,[\Mk/\Si\iin(\Mk)]=0\ $ 
holds in the Grothendieck group\/ $\Grot(\EndHMk)$.
\endproclaim

\Proof.
Consider the parabolic subalgebra $\calH\iin(z)_Q$ of $\calH_n(z)_Q$ generated 
by all elements $T_j$ with $0<j<n$, $j\ne i$. Let $Q\iin$ be the field $Q$ 
regarded as the 1-dimensional $\calH\iin(z)_Q$-module on which $T_j$ operates 
as multiplication by $z$ for $j<i$ and as $\,-\Id\,$ for $j>i$. By Lemma 2.3
$$
\sum_{i=0}^n\,(-1)^i\,[\Mk/\Si\iin(\Mk)]=
d\bigl(\,\sum_{i=0}^n\,(-1)^i\,[Q\iin\ot_{\calH_{i,n-i}(z)_Q}M_Q]\,\bigr).
$$
Since the algebra $\calH_n(z)_Q$ is semisimple, the functor 
$?\ot_{\calH_n(z)_Q}M_Q$ is exact, and therefore this functor induces a 
group homomorphism
$$
g:\ \Grotr\calH_n(z)_Q\to\Grot(\EndHMQ)
$$
where we denote by $\Grotr\calH_n(z)_Q$ the Grothendieck group of the category 
of finite dimensional right $\calH_n(z)_Q$-modules. Lemma 1.1 in its 
equivalent formulation for right modules shows that
$$
\sum_{i=0}^n\,(-1)^i\,[Q\iin\ot_{\calH_{i,n-i}(z)_Q}\calH_n(z)_Q]=0
\quad\hbox{in $\,\Grotr\calH_n(z)_Q$}.
$$
Applying the above map $g$, we get
$$
\sum_{i=0}^n\,(-1)^i\,[Q\iin\ot_{\calH_{i,n-i}(z)_Q}M_Q]=0
\quad\hbox{in $\,\Grot(\EndHMQ)$},
$$
and the desired conclusion follows.
\endproof

\proclaim
Lemma 2.5.
There are uniquely determined elements $\ze_\la\in\Grot(\EndHMk)$
with $\la\in\Par(n)$ such that
$$
[\,\bbk\triv\ot_{\calH_\mu(q)}\Mk\,]=\sum_{\la\in\Par(n)}K_{\la\mu}\,\ze_\la
\quad\hbox{for each $\,\mu\in\Par(n)$}.
\vadjust{\vskip-4pt} 
$$
Moreover{\rm,} $\ze_\la\ge0${\rm,} i.e. $\ze_\la$ represents an actual 
module{\rm,} for each $\la\in\Par(n)$.
\endproclaim

\Proof.
Let $S^\la$, $\la\in\Par(n)$, be the right Specht modules for the algebra 
$\calH_n(z)_Q$, as defined by Dipper and James \cite{Dip-J86}. By Lemma 1.1
$$
[Q\triv\ot_{\calH_\mu(z)_Q}\calH_n(z)_Q]=\sum_{\la\in\Par(n)}K_{\la\mu}\,[S^\la]
\quad\hbox{in $\Grotr\calH_n(z)_Q$}.
\vadjust{\vskip-4pt} 
$$
Applying the group homomorphism $g$ defined in the proof of Lemma 2.4, we get
$$
[Q\triv\ot_{\calH_\mu(z)_Q}M_Q]
=\sum_{\la\in\Par(n)}K_{\la\mu}\,[S^\la\ot_{\calH_n(z)_Q}M_Q]
\quad\hbox{in $\Grot(\EndHMQ)$}.
\vadjust{\vskip-4pt} 
$$
Applying now the decomposition map $d$ and making use of Lemma 2.2, we get the 
desired equalities in the group $\Grot(\EndHMk)$ with
$$
\ze_\la=d([S^\la\ot_{\calH_n(z)_Q}M_Q]).
$$
We have $\ze_\la\ge0$ since $d$ preserves positivity by the construction. 
Uniqueness of this collection of elements follows from the fact that the 
Kostka matrix is invertible. Indeed, this matrix is even unitriangular with 
respect to a suitable ordering of partitions \cite{Mac, Ch. I, (6.5)}.
\endproof

In the next lemma the assumption that the discrete valuation ring $O$ is 
complete goes into action.

\proclaim
Lemma 2.6.
If $X$ is any direct summand of the $\calH_n(q)$-module $\Mk,$ then 
$X=\Nk$ for some direct summand $N$ of the $\calH_n(z)$-module $M$.
\endproclaim

\Proof.
There exists an idempotent $e\in\EndHMk$ such that $X=\Im e$. 
By Lemma 2.1 $\EndHMk\cong\Ak\cong A/mA$ where 
$A=\EndHM$ and $m$ is the maximal ideal of $O$. The $O$-algebra $A$ 
is free of finite rank as a module. It is known that in this situation every 
idempotent of $A/mA$ lifts to an idempotent of $A$. Hence there exists an 
idempotent $\widetilde e\in A$ such that $e=\widetilde e\ot_O\bbk$. We may 
take $N=\Im\widetilde e$.
\endproof

\proclaim
Lemma 2.7.
Suppose that $X$ is a finite dimensional $\calH_n(q)$-module whose 
indecomposable direct summands all have $1$-dimensional sources. Then 
$X=\Mk$ for some $M\in\Rep$. If $q=-1,$ then we can even take $M\in\Triv$.  
\endproclaim

\Proof.
The class of $\calH_n(q)$-modules for which the conclusion of this lemma holds 
is obviously closed under direct sums. By Lemma 2.6 this class is closed also 
under direct summands. Since for $X=\calH_n(q)\ot_{\calH_{\la,\mu}(q)}\bbk_\lamu$ 
we can take $M=M^\lamu$, the first assertion holds then in full generality 
provided $q\ne-1$. If $q=-1$, then any $1$-dimensional representation of a 
parabolic subalgebra of $\calH_n(q)$ is trivial, so that we need only to look 
at the modules $X=\calH_n(q)\ot_{\calH_\la(q)}\bbk\triv$ for which $M=M^\la$ 
will do.
\endproof

\setitemsize(ii)
\proclaim
Corollary 2.8.
If $X$ is as in Lemma\/ {\rm2.7,} then{\rm:}

\item(i)
$\sum\limits_{i=0}^n\,(-1)^i\,[X/\Si\iin(X)]=0\ $ 
in the Grothendieck group $\,\Grot(\EndHX),$

\item(ii)
there exist finite dimensional\/ $\EndHX$-modules $V^\la,$ 
$\,\la\in\Par(n),$ such that in that group 
$\,[X/\Si_\mu(X)]=\sum_{\la\in\Par(n)}K_{\la\mu}\,[V^\la]\,$
for each $\mu\in\Par(n)$.

\endproclaim

This corollary repeats the conclusions of Lemmas 2.4 and 2.5. Note that 
$$
X/\Si_\mu(X)\cong\bbk\triv\ot_{\calH_\mu(q)}X.
$$

\section
3. The Hilbert series of the $R$-symmetric algebras

Let $R$ be a Hecke symmetry with parameter $q$ on a finite dimensional vector 
space $V$ over the ground field $\bbk$. For each $n\ge0$ consider the 
$\calH_n(q)$-module structure on $V^{\ot n}$ arising from $R$. Our first goal 
in this section is to describe a ring homomorphism $\ph:\Sym\to\Grot(R)$. All 
essential arguments needed to establish properties of $\ph$ are provided by 
the results of section 2.

As a preliminary step we will determine certain quotients of $V^{\ot n}$. In 
accordance with the notation introduced in section 2 we have for 
$(\la,\mu)\in\Par^2(n)$
$$ 
\Si_\lamu(V^{\ot n})=\sum_{i\in\calI^0_\lamu}\Im(R_i^{(n)}-q\cdot\Id)
+\sum_{i\in\calI^1_\lamu}\Ker(R_i^{(n)}-q\cdot\Id).
$$
This is an $\EndHn V^{\ot n}$-submodule, i.e. an $A_n(R)$-subcomodule, of 
$V^{\ot n}$. For each partition $\la=(\la_1,\ldots,\la_k)$ put
$$
\bbS^\la=\bbS_{\la_1}(V,R)\ot\ldots\ot\bbS_{\la_k}(V,R),\qquad
\La^\la=\La_{\la_1}(V,R)\ot\ldots\ot\La_{\la_k}(V,R).
$$
It will be assumed that $\bbS^\la=\La^\la=\bbk$ when $\la=0$.

\proclaim
Lemma 3.1.
For each $(\la,\mu)\in\Par^2(n)$ there is an isomorphism of $A_n(R)$-comodules 
$$
V^{\ot n}/\,\Si_\lamu(V^{\ot n})\cong\bbS^\la\ot\La^\mu.
$$
\endproclaim

\Proof.
Put $l=|\la|$ and $m=|\mu|$. Writing $V^{\ot n}$ as $V^{\ot l}\ot V^{\ot m}$ 
and noting that\hfil\break
$\Si_\lamu(V^{\ot n})=\Si_\lazero(V^{\ot l})\ot V^{\ot m}+
V^{\ot l}\ot\Si_\zeromu(V^{\ot m})$, we get
$$
V^{\ot n}/\,\Si_\lamu(V^{\ot n})\cong
V^{\ot l}/\,\Si_\lazero(V^{\ot l})\,\ot\,V^{\ot m}/\,\Si_\zeromu(V^{\ot m}).
$$
The ideal $I^\bbS$ defining the factor algebra $\bbS(V,R)$ of the tensor 
algebra $\bbT(V)$ has homogeneous components 
$I^\bbS_k=\sum_{0<i<k}\Im(R_i^{(k)}-q\cdot\Id)\sbs V^{\ot k}$. Hence
$$
\Si_\lazero(V^{\ot l})=\sum_{j=1}^{\ell(\la)}V^{\ot(\la_1+\ldots+\la_{j-1})}\ot 
I^\bbS_{\la_j}\ot V^{\ot(\la_{j+1}+\ldots+\la_{\ell(\la)})},
$$
and it follows that $V^{\ot l}/\,\Si_\lazero(V^{\ot l})\cong\bbS^\la$.
On the other hand, the ideal $I^\La$ defining $\La(V,R)$ has homogeneous 
components $I^\La_k=\sum_{0<i<k}\Ker(R_i^{(k)}-q\cdot\Id)$. One obtains 
similarly $V^{\ot m}/\,\Si_\zeromu(V^{\ot m})\cong\La^\mu$.
\endproof

\proclaim
Proposition 3.2.
Suppose that $R$ satisfies the $1$-dimensional source condition. Then there is 
a ring homomorphism $\ph:\Sym\to\Grot(R)$ such that

\item(i)
$\,\ph(h_n)=[\bbS_n(V,R)]$ and $\ph(e_n)=[\La_n(V,R)]$ for all $n\ge0,$

\item(ii)
$\,\ph(s_\la)\ge0$ for all $\la\in\Par$.

\endproclaim

\Proof.
Since $h_1,h_2,\ldots$ are algebraically independent generators of the ring 
of symmetric functions $\Sym$ (see \cite{Mac, (2.8)}), homomorphisms from $\Sym$ 
to another ring are uniquely determined by their values on those elements. 
Thus we can define $\ph$ setting $\ph(h_n)=[\bbS_n(V,R)]$ for all $n>0$. Since 
$h_0=e_0=1$ and $\bbS_0(V,R)=\La_0(V,R)=\bbk$ is the trivial 1-dimensional 
$A(R)$-comodule 
which represents the identity element of the ring $\Grot(R)$, the 
two equalities in (i) are obvious for $n=0$. If $n>0$, then the relation 
$\sum_{i=0}^n(-1)^ih_ie_{n-i}=0$ in $\Sym$ yields 
$\sum_{i=0}^n(-1)^i\ph(h_i)\ph(e_{n-i})=0$. The value $\ph(e_n)$ 
is determined by induction on $n$ from the latter equality. To show that the 
second equality in (i) holds for all $n>0$ we have to check that
$$
\sum_{i=0}^n\,(-1)^i[\mskip2mu\bbS_i(V,R)]\cdot[\La_{n-i}(V,R)]
=\sum_{i=0}^n\,(-1)^i[\mskip2mu\bbS_i(V,R)\ot\La_{n-i}(V,R)]=0.\eqno(\rm i')
$$

If $\mu\in\Par(n)$, then $\ph(h_\mu)=\Smu$. The equality 
$h_\mu=\sum_{\la\in\Par(n)}K_{\la\mu}s_\la$ in $\Sym$ entails 
$[\mskip2mu\Smu]=\sum_{\la\in\Par(n)}K_{\la\mu}\ph(s_\la)$. This determines 
the values of $\ph$ on the Schur functions since the Kostka matrices are 
invertible. Part (ii) means that there exists a collection of finite dimensional 
$A(R)$-comodules $V^\la$, $\,\la\in\Par$, such that $\ph(s_\la)=[V^\la]$ for 
each $\la\in\Par$. Those equalities are equivalent to the equalities
$$
[\mskip2mu\Smu]=\sum_{\la\in\Par(n)}K_{\la\mu}\,[V^\la]\qquad\hbox{for each 
$\mu\in\Par(n)$}.\eqno(\rm ii')
$$
The validity of both (i$'$) and (ii$'$) in the group $\Grot_n(R)$ is ensured 
by Corollary 2.8. Indeed, for $X=V^{\ot n}$ the group $\Grot\EndHX$ has been 
identified with $\Grot_n(R)$, and we have
$$
X/\Si\iin(X)\cong\bbS_i(V,R)\ot\La_{n-i}(V,R)\quad{\rm and}\quad
X/\Si_\mu(X)\cong\Smu
$$
by Lemma 3.1.
\endproof

\proclaim
Corollary 3.3.
In the ring $\,\Grot(R)[[t]]$ consider the formal power series
$$
G_\bbS=\sum_{i=0}^\infty\,\,[\mskip2mu\bbS_i(V,R)]\,t^i,\qquad
G_\La=\sum_{i=0}^\infty\,\,[\La_i(V,R)]\,t^i.
$$
We have $\,G_\bbS(t)G_\La(-t)=1$.
\endproclaim

\Proof.
The coefficient of $t^n$ in the product $G_\bbS(t)G_\La(-t)$ vanishes for 
each $n>0$ according to equality (i$'$) in the proof of Proposition 3.2. The 
constant term of this product is the identity element 
$[\mskip1mu\bbk\mskip1mu]$ of the ring $\Grot(R)$.
\endproof

Corollary 3.3 strengthens the well known relation $\,H_\bbS(t)H_\La(-t)=1\,$ 
between the Hilbert series of the algebras $\bbS=\bbS(V,R)$ and $\La=\La(V,R)$. 
This relation was proved by Gurevich \cite{Gur90} in the semisimple case. The 
more general case when $q$ is arbitrary and $R$ satisfies the $1$-dimensional 
source condition was treated in \cite{Sk19}.

\medbreak
For each $\la\in\Par$ we will denote by $V^\la$ any finite dimensional right 
$A(R)$-comodule such that $[V^\la]=\ph(s_\la)$. In the nonsemisimple case such 
a comodule is generally not unique, but we will be concerned with numeric 
characteristics such as the dimension of $V^\la$ which depend only on the 
class of $V^\la$ in the group $\Grot(R)$.

Let $K$ be a commutative ring. For a formal power series 
$f=\sum_{i=0}^\infty a_it^i\in K[[t]]$ with $a_0=1$ we denote by 
$$
\fhat:\Sym\mapr{}K
$$
the ring homomorphism such that $h_n\mapsto a_n$ for each $n\ge0$. We will 
assume tacitly that $a_i=0$ for all integers $i<0$.

\proclaim
Corollary 3.4.
Let $f=\sum_{i=0}^\infty a_it^i\in\bbZ[[t]]$ be the Hilbert series of the 
$R$-symmetric algebra $\bbS(V,R)$. Then
$$
\dim V^\la=\fhat(s_\la)=\det\,(a_{\la_i-i+j})_{1\le i,j\le k}
$$
for each partition $\la$ with $k$ parts $\la_1,\ldots,\la_k$.
\endproclaim

\Proof.
There is a ring homomorphism $\de:\Grot(R)\to\bbZ$ such that $\de([X])=\dim X$ 
for each finite dimensional right $A(R)$-comodule $X$. We have\/ 
$\fhat=\de\circ\ph$ since the left and right hand sides of this equality are 
presented by ring homomorphisms with the same values $a_n$ on the 
generators $h_n$ of the ring $\Sym$. 

Hence $\fhat(s_\la)=\de\bigl(\ph(s_\la)\bigr)=\dim V^\la$. 
The second equality follows from the Jacobi-Trudi identity 
$\,s_\la=\det\,(h_{\la_i-i+j})_{1\le i,j\le k}\,$.
\endproof

An infinite sequence of real numbers $a_0,a_1,a_2,\ldots$ with $a_0=1$ is called 
\emph{totally positive}\/ or a \emph{P\'olya frequency sequence}\/ if all minors 
of the infinite Toeplitz matrix $(a_{j-i})_{i,j\ge0}$ are nonnegative. Denoting 
by $f$ the generating series of the given sequence, we see that $\fhat(s_\la)$ 
is one of these minors taken from a set of consecutive columns. Other minors 
are the values of $\fhat$ on skew Schur functions $s_{\nu/\mu}$. Since 
$s_{\nu/\mu}=\sum_{\la\in\Par}c_{\la\mu}^\nu s_\la$ with the 
Littlewood-Richardson coefficients $c_{\la\mu}^\nu\ge0$, it 
follows that the sequence is totally positive if and only if
$$
\fhat(s_\la)\ge0\quad\hbox{for each $\la\in\Par$},
$$
i.e. there is no need to look at the other minors. This formulation of total 
positivity together with two other equivalent conditions is discussed by 
Stanley \cite{Stan, Exercise 7.91e} in the case when $f$ is a polynomial. The 
statement given above is equally valid for infinite series.

Since $\dim V^\la\ge0$ for all $\la\in\Par$, it follows from Corollary 3.4 
that the dimensions of the homogeneous components $\bbS_n(V,R)$ form a 
totally positive sequence. All possibilities for the generating series $f$ of 
a totally positive sequence were determined in the work of Aissen, Schoenberg, 
Whitney \cite{Ais-SW52} and Edrei \cite{Edr52}. In particular, \cite{Ais-SW52, 
Th.~1} states that $f$ converges in a neighborhood of $0$ in $\bbC$ and 
extends to a meromorphic function with negative real zeros and positive real 
poles on the whole $\bbC$. As it turns out, $f$ is rational precisely when 
$\fhat(s_\la)=0$ for at least one $\la\in\Par$.

Thus rationality of the Hilbert series of $\bbS(V,R)$ is equivalent to the 
existence of partitions $\la$ for which $V^\la=0$. As observed by \PHHai 
\cite{Hai99}, in the semisimple case such $\la$ do exist since otherwise the 
representation of $\calH_n$ in $V^{\ot n}$ would be faithful for each $n$, 
but this is impossible since $\dim\calH_n=n!>(\dim V)^{2n}$ for large $n$. 
We cannot use this argument directly when $q$ is a root of 1. However, the 
next lemma provides a replacement.

For integers $a,k>0$ denote by $(a^k)$ the partition $(a,\ldots,a)$ with 
exactly $k$ parts, each equal to $a$.

\proclaim
Lemma 3.5.
Suppose that $R$ satisfies the $1$-dimensional source condition. There exist  
integers $n>0$ and $k>0$ such that $V^{(n^k)}=0$. Moreover{\rm,} if 
$V^{(n^k)}=0,$ then $V^{(m^k)}=0$ for all $m>n$.
\endproclaim

\Proof.
Denote by $\dla$ the dimension of the Specht module $S^\la$, i.e. $\dla$ is 
equal to the Kostka number $K_{\la,(1^n)}$ counting the standard $\la$-tableaux. 
We first note that $V^\la=0$ for each $\la\in\Par(n)$ with $\dla>(\dim V)^n$. 
To prove this we start with the equality
$$
h_1^n=h_{(1^n)}=\sum_{\rho\mskip1mu\in\Par(n)}\drho s_\rho\quad
\hbox{in the ring $\,\Sym$}.
$$
An application of $\ph$ yields $\,[V^{\ot n}]
=\sum_{\rho\mskip1mu\in\Par(n)}\drho\,[V^\rho]\,$ in $\Grot(R)$. Hence
$$
(\dim V)^n=\dim V^{\ot n}=\sum_{\rho\mskip1mu\in\Par(n)}\drho\dim V^\rho,
$$
and therefore $\dla\dim V^\la\le(\dim V)^n$, which entails the previous claim.

By the hook length formula
$$
d^{(n^k)}
=\frac{(kn)!}{\prod\limits_{i=1}^k\prod\limits_{j=1}^n(n+k-i-j+1)}
=(kn)!\prod_{i=1}^k\frac{(k-i)!}{(n+k-i)!}
=\frac{(kn)!}{n!^k}\prod_{i=0}^{k-1}{n+i\choose i}^{-1}\!.
$$
Using the Stirling asymptotic formula $n!\sim n^ne^{-n}\sqrt{2\pi n}$ and 
observing that ${n+i\choose i}\sim n^i/i!$ for each $i$, we deduce that
$$
d^{(n^k)}\sim k^{kn}n^{-(k^2-1)/2}\cdot\sqrt{k/(2\pi)^{k-1}}\,\prod_{i=0}^{k-1}i!
\qquad\hbox{as $n\to\infty$}.
$$
If $k>\dim V$, then $d^{(n^k)}>(\dim V)^{kn}$ for large $n$. As we have 
observed, this yields $V^{(n^k)}=0$, proving the first assertion of Lemma 3.5.

Next, if $\la\in\Par$ is such that $V^\la=0$, then 
$0=[V^\la][V^\mu]=\sum_{\nu\in\Par}c_{\la\mu}^\nu[V^\nu]$ for each $\mu\in\Par$. Since 
$c_{\la\mu}^\nu\ge0$ for all $\nu$, it follows that $V^\nu=0$ whenever 
$c_{\la\mu}^\nu\ne0$. In particular, this applies in the case when 
$\la=(n^k)$, $\mu=((m-n)^k)$ and $\nu=\nobreak(m^k)$ with $m>n$. Indeed, for 
these partitions we have $c_{\la\mu}^\nu=1$ by the Littlewood-Richardson rule.  
\endproof

For the application to the Hilbert series of the $R$-symmetric algebras the 
full strength of the analytic result on totally positive sequences is actually 
not needed. For one thing rationality of the series follows from a purely 
algebraic fact formulated in terms of Hankel determinants (see \cite{Hen, 
Th. 7.5f}). These determinants with the reversed order of rows are certain 
minors of the Toeplitz matrix. The next lemma gives a version of that result 
paying attention to the ring of coefficients.

\proclaim
Lemma 3.6.
Let $f=\sum_{i=0}^\infty a_it^i$ be a formal power series with coefficients in 
a commutative Noetherian domain $K$. For each pair of integers $i\ge k>0$ put
$$
\De_i^{(k)}=\fhat(s_{(i^k)})
=\left|\matrix{a_i & a_{i+1} & \cdots & a_{i+k-1} \cr
a_{i-1} & a_i & \cdots & a_{i+k-2} \cr
\vdots & \vdots & & \vdots \cr
a_{i-k+1} & a_{i-k+2} & \cdots & a_i \cr
}\right|.
$$
Suppose that there are integers $n>r>0$  such that $\De_i^{(r)}\ne0$ and 
$\De_i^{(r+1)}=0$ for all $i\ge n$. Then $f=q^{-1}p$ for some polynomials 
$p,q\in K[t]$ with $q(0)=1$. If $K$ is integrally closed{\rm,} then such an 
expression holds with $q$ of degree $r,$ in which case $p$ and $q$ are 
relatively prime in the ring $Q[t]$ where $Q$ is the field of fractions of $K$.
\endproclaim

\Proof.
Put $A_i^{(k)}=(a_{i-k+1},a_{i-k+2}\,,\ldots,a_i)$. Part of the hypothesis 
means that for each $i\ge n$ the $r$ vectors 
$A_i^{(r)}\!,A_{i+1}^{(r)},\ldots,A_{i+r-1}^{(r)}$ are linearly independent 
over $Q$ and so form a basis for the $r$-dimensional vector space $Q^r$. Then 
$$
A_{i+r}^{(r)}=c_1A_{i+r-1}^{(r)}+\ldots+c_{r-1}A_{i+1}^{(r)}+c_rA_i^{(r)}
\eqno({\rm Rel}_i)
$$
with uniquely determined coefficients $c_1,\ldots,c_r\in Q$. We claim that 
these coefficients do not depend on $i$. This will follow once we show that 
for each $i>n$ the coefficients in $({\rm Rel}_i)$ are the same as those in 
$({\rm Rel}_{i-1})$. But $A_{i+r}^{(r+1)}$ is a linear combination of vectors 
$A_i^{(r+1)},A_{i+1}^{(r+1)},\ldots,A_{i+r-1}^{(r+1)}\in Q^{r+1}$ because 
$\De_i^{(r+1)}=0$. Since the projection $Q^{r+1}\to Q^r$ onto the last $r$ 
components maps $A_j^{(r+1)}$ to $A_j^{(r)}$ for each $j$, we must have
$$
A_{i+r}^{(r+1)}=c_1A_{i+r-1}^{(r+1)}+\ldots+c_{r-1}A_{i+1}^{(r+1)}+c_rA_i^{(r+1)}
$$
with coefficients from $({\rm Rel}_i)$. Noting that the projection 
$Q^{r+1}\to Q^r$ onto the first $r$ components maps $A_j^{(r+1)}$ to 
$A_{j-1}^{(r)}$ for each $j$, we get $({\rm Rel}_{i-1})$ with the same 
coefficients, and the claim is proved.

Now $({\rm Rel}_i)$ shows that $a_{i+1}=\sum_{j=1}^rc_ja_{i+1-j}$ for each 
$i\ge n$. Setting 
$$
g=1-{\textstyle\sum\limits_{j=1}^r}c_jt^j\in Q[t],
$$
we see that the coefficient of $t^{i+1}$ in the formal power series $gf$ 
vanishes whenever $i\ge n$. Hence $gf$ is a polynomial with coefficients in $Q$.

Define a linear operator $\th$ on the vector space $Q^r$ by the formula
$$
\th(x_1,x_2,\ldots,x_r)=(x'_1,x'_2,\ldots,x'_r)
$$
where $x'_j=x_{j+1}$ for $0<j<r$ and $x'_r=\sum_{j=1}^rc_jx_{r+1-j}$. We have 
$\th(A_i^{(r)})=A_{i+1}^{(r)}$ for each $i\ge n$, and it follows that
$$
\th^r-c_1\th^{r-1}-\ldots-c_{r-1}\th-c_r\Id=0
$$
since this operator annihilates all vectors $A_i^{(r)}$ with $i\ge n$ in view 
of $({\rm Rel}_i)$. On the other hand, the operators 
$\th,\ldots,\th^{r-1},\th^r$ are linearly independent over $Q$ since so are 
the vectors $A_{i+1}^{(r)},\ldots,A_{i+r-1}^{(r)},A_{i+r}^{(r)}$ for $i\ge n$.  
Hence $c_r\ne0$, and $g$ is a scalar multiple of the minimal polynomial of the 
inverse operator $\th^{-1}$.

Note that $\th(M)\sbs M$ where $M$ is the submodule of the free $K$-module 
$K^r$ generated by $\{A_i^{(r)}\mid i\ge n\}$. Since $K$ is Noetherian, $M$ 
has to be finitely generated. Hence $\th$ is integral over $K$, i.e. for some 
$k>0$ there exists a relation
$$
\th^k+e_1\th^{k-1}+\ldots+e_{k-1}\th+e_k\Id=0
$$
with coefficients $e_j\in K$. Take $q=1+\sum_{j=1}^ke_jt^k$. Then $g$ divides 
$q$ since $q(\th^{-1})=0$, and it follows that $qf$ is a polynomial whose 
coefficients are in $K$ since so are the coefficients of $f$ and $q$. With 
$p=qf$ the first conclusion is thus proved.

We can write $q=(1-\xi_1t)\cdots(1-\xi_kt)$ with $\xi_1,\ldots,\xi_k$ in the 
algebraic closure of the field $Q$. Each $\xi_j$ is integral over $K$ since 
$q(\xi_j^{-1})=0$. Since $g$ is a divisor of $q$ with $g(0)=1$, it is the 
product of some of these factors $1-\xi_jt$. Hence all coefficients of $g$ are 
integral over $K$ too. If $K$ is integrally closed, we get $g\in K[t]$, and so 
we may take $q=g$. Then $p$ and $q$ cannot have a common divisor in $Q[t]$ 
since otherwise $hf\in Q[t]$ for some polynomial $h\in Q[t]$ of degree less 
than $r$, but this implies that the sequence $(a_i)$ starting at some term 
satisfies a linear recurrence relation of order less than $r$, which 
contradicts the linear independence of the previously considered vectors 
$A_i^{(r)}\!,A_{i+1}^{(r)},\ldots,A_{i+r-1}^{(r)}$.
\endproof

Knowing rationality of $f$, one needs only part of the arguments given in 
\cite{Ais-SW52} to determine the location of zeros and poles. Moreover, 
Stanley's criterion of total positivity provides further simplifications:

\proclaim
Lemma 3.7.
Let $f=\sum_{i=0}^\infty a_it^i\in\bbR[[t]]$ be the generating series of a 
totally positive sequence. Suppose that $f$ represents a rational function of 
$t$. Then all its zeroes are negative{\rm,} while all its poles are positive 
real numbers.
\endproclaim

\Proof.
If $a_n\ne0$ and $a_{n+1}=0$ for some $n\ge0$, then $f$ is a polynomial since 
for each $i>n$ the conditions $a_ia_{n+1}-a_{i+1}a_n\ge0$ and $a_{i+1}a_n\ge0$ 
entail $a_{i+1}=0$.

Suppose that $a_n\ne0$ for all $n\ge0$. Then the positive numbers $a_{n+1}/a_n$ 
form a monotone nonincreasing sequence. Hence this sequence converges to some 
$\ga\ge0$. Rationality of $f$ implies that $\ga^{-1}$ is one of its poles, i.e. 
$1-\ga t$ is a divisor of the denominator $q$ in the expression of $f$ as a 
fraction of two relatively prime polynomials. In particular, $\ga\ne0$. The 
power series $f_\flat=(1-\ga t)f$ has coefficients
$$
a'_i=a_i-\ga a_{i-1}
$$
which form a totally positive sequence. Indeed, given a partition 
$\la=(\la_1,\ldots,\la_k)$, we have
$$
\fhat_\flat(s_\la)=\det\,(a'_{\la_i-i+j})_{1\le i,j\le k}
=\left|\matrix{1&\ga&\cdots&\ga^k\cr
a_{\la_1-1} & a_{\la_1} & \cdots & a_{\la_1+k-1}\cr
\vdots & \vdots & & \vdots \cr
a_{\la_k-k+1} & a_{\la_k-k+2} & \cdots & a_{\la_k}\cr}\right|.
$$
Since $\,\ga^i=\lim_{n\to\infty}a_{n+i}/a_n\,$, this yields 
$$
\fhat_\flat(s_\la)=\lim_{n\to\infty}\fhat(s_{\la^{(n)}})/a_n\ge0
$$
where we put $\la^{(n)}=(n,\la_1,\ldots,\la_k)$. Thus $f_\flat$ satisfies the 
same assumptions as $f$. We have seen that $f$ has a pole $\ga^{-1}>0$. 
Proceeding by induction on the degree of the denominator $q$, we conclude that 
all poles of $f$ are positive.

The coefficients $b_i$ of the power series $g(t)=1/f(-t)$ also form a totally 
positive sequence. This fact was proved in \cite{Ais-SW52}, but again it can 
be explained very easily within the theory of symmetric functions. Since 
$\sum_{i=0}^n(-1)^ia_ib_{n-i}=0$ for $n>0$, we have $b_n=\fhat(e_n)$. Hence 
$\ghat=\fhat\circ\om$ where $\om$ is the automorphism of the ring $\Sym$ such 
that $h_n\mapsto e_n$ for each $n$. Since $\om(s_\la)=s_{\la'}$ where $\la'$ 
is the conjugate of $\la$, we get $\ghat(s_\la)=\fhat(s_{\la'})\ge0$ for each 
$\la\in\Par$. Thus all poles of $g$ are positive, which means that all zeroes 
of $f$ are negative.
\endproof

\proclaim
Theorem 3.8.
Suppose that $R$ satisfies the $1$-dimensional source condition. Then
$$
H_{\La(V,R)}(t)=f_0(-t)/f_1(t),\qquad H_{\bbS(V,R)}(t)=f_1(-t)/f_0(t)
$$
with integer polynomials $\,f_0,\,f_1\in\bbZ[t]\,$ whose constant terms are 
equal to\/ $1$ and all roots are positive real numbers.
\endproclaim

\Proof.
Put $f=H_{\bbS(V,R)}$. Since $H_{\La(V,R)}(t)=f(-t)^{-1}$, it suffices to prove 
only the formula for $f$. Define $\De_i^{(k)}$ as in Lemma 3.6. By Corollary 
3.4 $f$ is the generating series of a totally positive sequence, and also 
$\De_i^{(k)}=\dim V^{(i^k)}$. Let $r$ be the smallest nonnegative integer for 
which there exists $n>0$ such that $V^{(n^{r+1})}=0$. By Lemma 3.5 $r$ is 
well-defined and $V^{(m^{r+1})}=0$ for all $m>n$.

If $r=0$, then $\bbS_n(V,R)=0$, which means that $f$ is a polynomial. Note that 
the constant term of $f$ is equal to 1.

If $r>0$, then the determinants $\De_i^{(k)}$ satisfy the assumption of Lemma 
3.6. Taking $K=\bbZ$, we deduce that $f$ is a fraction of two relatively prime 
integer polynomials with constant terms equal to 1. By Lemma 3.7 $f$ has 
negative zeros and positive poles. This ensures the desired properties of $f_0$ 
and $f_1$.
\endproof

Now we will extend to the present situation two additional results obtained by 
\PHHai in the semisimple case \cite{Hai99, Th. 5.1, Cor. 5.2}.

\proclaim
Corollary 3.9.
Let $(r_0,r_1)$ be the birank of $R,$ i.e. $r_i=\deg f_i$ for $i=0,1$. 
Then $r_0+r_1\le\dim V$. Moreover{\rm,} if $r_0+r_1=\dim V,$ then
$$
H_{\La(V,R)}(t)=(1+t)^{r_0}/(1-t)^{r_1},\qquad 
H_{\bbS(V,R)}(t)=(1+t)^{r_1}/(1-t)^{r_0}
$$
\endproclaim

\Proof.
We have $f_0(t)=\prod_{i=1}^{r_0}(1-\al_it)$ and 
$f_1(t)=\prod_{i=1}^{r_1}(1-\be_it)$ with $\al_i,\be_i>0$. The fact that all 
coefficients of these polynomials are integers entails
$$
{\al_1+\ldots+\al_{r_0}\over r_0}\ge
\root\raise2pt\hbox{$\scriptstyle r_0$}\of{\al_1\cdots\al_{r_0}\strut}\ge1,
\qquad{\be_1+\ldots+\be_{r_1}\over r_1}\ge
\root\raise2pt\hbox{$\scriptstyle r_1$}\of{\be_1\cdots\be_{r_1}\strut}\ge1.
$$
Since the coefficient of $t$ in $H_{\bbS(V,R)}$ is equal to the dimension of 
$V$, we get
$$
\dim V=\tsum\al_i+\tsum\be_j\ge r_0+r_1\,.
$$
If the equality is attained here, then $\sum\al_i=r_0$ and $\sum\be_j=r_1$, so 
that the equalities are attained also in the previously displayed formulas. 
This is only possible when all $\al_i$ and $\be_j$ are equal to 1.
\endproof

\proclaim
Corollary 3.10.
Let $(r_0,r_1)$ be the birank of $R$. Then $V^\la\ne0$ if and only if 
$\la_j\le r_1$ for all $j>r_0,$ i.e. $\la\in\Ga(r_0,r_1)$.
\endproclaim

\Proof.
Let $x_1,\ldots,x_{r_0},y_1,\ldots,y_{r_1}$ be commuting indeterminates. Consider 
the ring homomorphism $\Sym\to\bbZ[x_1,\ldots,x_{r_0},y_1,\ldots,y_{r_1}]$ under 
which the formal power series $\sum_{n=0}^\infty e_nt^n$ and 
$\sum_{n=0}^\infty h_nt^n$ specialize, respectively, to
$$
\prod_{i=1}^{r_0}(1+x_it)\cdot\prod_{j=1}^{r_1}(1-y_jt)^{-1}\quad{\rm and}\quad
\prod_{j=1}^{r_1}(1+y_jt)\cdot\prod_{i=1}^{r_0}(1-x_it)^{-1}.
$$ 
Denote by $u(x/y)$ the image of $u\in\Sym$ in the ring 
$\bbZ[x_1,\ldots,x_{r_0},y_1,\ldots,y_{r_1}]$ and by $u(\al/\be)$ the value of 
this polynomial $u(x/y)$ at the point 
$$
(\al_1,\ldots,\al_{r_0},\be_1,\ldots,\be_{r_1})\in\bbR^{r_0+r_1}
$$
where $\al_i$ and $\be_j$ are as in the proof of Corollary 3.9. Then
$$
\dim V^\la=s_\la(\al/\be)\quad\hbox{\rm for each $\la\in\Par$}.
$$
Since the left and right hand sides of this equality are evaluations at 
$s_\la$ of two ring homomorphisms $\Sym\to\bbZ$, it suffices to check it on 
the generators $s_{(n)}=h_n$ of the ring $\Sym$. But for $\la=(n)$ we have 
$V^\la=\bbS_n(V,R)$, while $s_\la(\al/\be)=h_n(\al/\be)$ is exactly the 
coefficient of $t^n$ in the Hilbert series of the algebra 
$\bbS(V,R)$, so that the equality is indeed true.

The specialization $u\mapsto u(x/y)$ defined above differs from that of Macdonald 
\cite{Mac, Ch.~I, Example 3.23} in that each $y_j$ is changed to $-y_j$. With 
this change formula (1) in \cite{Mac, Ch.~I, Example 5.23} reads as
$$
s_\la(x/y)=\sum_{\mu\in\Par}s_\mu(x)s_{\la'\!/\mu'}(y)
$$
where $s_{\la'\!/\mu'}$ is the skew Schur function corresponding to the pair 
$\la',\mu'$ of partitions conjugate to $\la$ and $\mu$. Thus $s_\la(x/y)$ is 
precisely the hook Schur function $HS_\la$ in the notation and terminology of 
Berele and Regev \cite{Ber-R87, Definition 6.3}.

Recall that each skew Schur function is a linear combination of monomial 
symmetric functions with nonnegative integer coefficients. Hence any monomial 
in the indeterminates $x_1,\ldots,x_{r_0},y_1,\ldots,y_{r_1}$ has 
nonnegative coefficient in $s_\la(x/y)$. It is known also that $s_\la(x/y)\ne0$ 
if and only if $\la\in\Ga(r_0,r_1)$ \cite{Ber-R87, Cor. 6.5}. Since all real 
numbers $\al_i$ and $\be_j$ are positive, we deduce that 
$\,\dim V^\la=s_\la(\al/\be)>0\,$ if and only if $\la\in\Ga(r_0,r_1)$.
\endproof

\Remark.
By Corollary 3.10 the image of $\ph$ is the subgroup of $\Grot(R)$ generated 
by the classes $[V^\la]$ with $\la\in\Ga(r_0,r_1)$. It is not clear whether 
these classes are always linearly independent over $\bbZ$. In any event for 
each $n\ge0$ the rank of the free abelian group $\ph(\Sym_n)$ does not exceed 
the cardinality of the set $\Ga(r_0,r_1)\cap\Par(n)$. On the other hand, the 
group $\Grot_n(R)$ may have a larger rank. When this happens, $\ph$ is not 
surjective unlike what we have seen in the semisimple case.

As an example consider the supersymmetry $R$ on a $\bbZ/2\bbZ$-graded vector 
space $V=V_0\oplus V_1$. It is defined by the rule $R(v\ot w)=w\ot v$ for 
homogeneous elements $v,w\in V$ at least one of which is even, and 
$R(v\ot w)=-w\ot v$ when both $v$ and $w$ are odd. We assume here that 
$\chr\bbk\ne2$. This operator is a Hecke symmetry with parameter $q=1$, so 
that $\calH_n(q)$ is the group algebra $\bbk\frS_n$. It is easy to see that 
the $1$-dimensional source condition is satisfied. Since $\bbS(V,R)$ is the 
tensor product $S(V_0)\ot\bigwedge(V_1)$ of the ordinary symmetric and 
exterior algebras, it has Hilbert series $(1+t)^{r_1}/(1-t)^{r_0}$ where 
$r_i=\dim V_i$ for $i=0,1$. Thus the birank of $R$ coincides with the 
superdimension of $V$.

The right $A_n(R)$-comodules may be identified with the degree $n$ polynomial 
representations of the supergroup $GL_{r_0\mid r_1}$ over the field $\bbk$. 
(In the framework of super theory one endows $A(R)=\bigoplus_{k=0}^\infty 
A_k(R)$ with a modified version of the multiplication described in section 1. 
This modified multiplication makes $A(R)$ into a super bialgebra rather than 
an ordinary bialgebra, and in this way $A(R)$ is identified with the 
subalgebra of the coordinate algebra of $GL_{r_0\mid r_1}$ generated by the 
coefficient functions of the natural representation on $V$.)

Thus the rank of $\Grot_n(R)$ equals the number of irreducible polynomial 
representations of $GL_{r_0\mid r_1}$ of degree $n$. Irreducible 
representations are classified by their highest weights with respect to a 
maximal torus $T$ of the group $GL_{r_0}\times GL_{r_1}\sbs GL_{r_0\mid r_1}$. 
Highest weights  of polynomial representations may be interpreted as pairs of 
partitions $\la,\mu$ such that $\ell(\la)\le r_0$ and $\ell(\mu)\le r_1$. 
However, not all such pairs correspond to a polynomial representation.

If $\chr\bbk=0$, then the highest weights of irreducible polynomial
representations of $GL_{r_0\mid r_1}$ are selected by the additional condition 
$\ell(\mu)\le\la_{r_0}$ on the pair $(\la,\mu)$ (see \cite{Ser84, Cor.~1 to 
Th.~2}). On the other hand, since the group algebras $\bbk\frS_n$ are semisimple, 
the irreducible polynomial $GL_{r_0\mid r_1}$-modules of degree $n$ are 
precisely the simple $A_n(R)$-comodules
$$
V^\nu_R=\Hom_{\bbk\frS_n}(S^\nu,V^{\ot n})
$$
with $\nu\in\Ga(r_0,r_1)\cap\Par(n)$, while $V^\nu_R=0$ when 
$\nu\notin\Ga(r_0,r_1)$. In terms of representations of Lie superalgebras this 
fact was established long ago independently by Sergeev \cite{Ser84} and 
Berele, Regev \cite{Ber-R87}. It should be noted also that the character of 
$V^\nu_R$ defined in terms of weight spaces with respect to the torus $T$ is 
exactly the hook Schur function $s_\nu(x/y)$. This explains the properties of 
these functions referred to in the proof of Corollary 3.10.

Suppose now that $\chr\bbk=p>0$. In this case the irreducible representation 
of $GL_{r_0\mid r_1}$ with highest weight represented by the pair $(\la,\mu)$ 
is polynomial if and only if $j(\mu)\le\la_{r_0}$ where $j(\mu)$ is the 
cardinality of a combinatorially defined subset of nodes of the Young diagram 
of $\mu$ which contains at most one node from each row of the diagram. This 
condition was found by Brundan and Kujawa \cite{Br-K03, Th. 6.5}. The 
inequality $j(\mu)\le\ell(\mu)$ holds for each $\mu$. At the same time there 
exist highest weights $(\la,\mu)$ such that $\ell(\mu)>\la_{r_0}$ but 
$j(\mu)\le\la_{r_0}$. For example, $j(\mu)=0$ if $\mu_i\equiv0\mod p$ for all 
$i$. From this it follows that the group $\Grot_n(R)$ has larger rank than 
that in the case of a field of characteristic 0, and therefore 
$\ph(\Sym_n)\ne\Grot_n(R)$, for infinitely many $n$.

\endremark

\section
4. Tensor powers $V^{\ot n}$ as modules over the Hecke algebras

Consider $V^{\ot n}$ as an $\calH_n(q)$-module with respect to the 
representation arising from a Hecke symmetry $R$. Assuming that $R$ satisfies 
the $1$-dimensional source condition, we will associate with this module a 
symmetric function $\ch([V^{\ot n}])\in\Sym_n$ which encodes enough 
information to determine the image of $V^{\ot n}$ in the Grothendieck group 
$\Grot\calH_n(q)$.

For an associative algebra $\frA$ over some field, say $F$, denote by 
$\,\RepA\,$ the abelian group generated by the isomorphism classes $[X]$ of 
finite dimensional left $\frA$-modules with the defining relations 
$[X]=[X']+[X'']$ for each triple of finite dimensional left $\frA$-modules 
such that $X\cong X'\oplus X''$. It is a free abelian group with a basis 
consisting of the isomorphism classes of indecomposable finite dimensional 
left $\frA$-modules. By abuse of notation $[X]$ will stand for an element of 
either $\RepA$ or $\Grot\frA$, depending on the context. There is a canonical 
group homomorphism
$$
c:\RepA\mapr{}\Grot\frA
$$
sending the class of $X$ in $\RepA$ to the class of $X$ in $\Grot\frA$. This 
map is an isomorphism when $\frA$ is semisimple. In general $X\cong Y$ 
whenever $[X]=[Y]$ in $\RepA$.

Consider a $\bbZ$-bilinear form on $\RepA$ defined by the formula
$$
\<\,[X],[Y]\,\>=\<X,Y\>=\dim_F\Hom_{\mskip1mu\frA}(X,Y)
$$
for each pair of finite dimensional left $\frA$-modules $X$ and $Y$.

Denote by $\RepHq$ (respectively, $\TrivHq$) the subgroup of $\RHq$ generated 
by the isomorphism classes of indecomposable left $\calH_n(q)$-modules which 
have a 1-dimensional (respectively, trivial) source. Then 
$\,\TrivHq\sbs\RepHq$.

In the next lemma we work in the settings of section 2:

\medbreak
\proclaim
Lemma 4.1.
There is a group homomorphism $\,e:\RepHq\to\Grot\calH_n(z)_Q\,$ such that 
$\,e([\Mk])=[M_Q]\,$ for each lattice $M\in\Rep$ if $q\ne-1$ and for 
$M\in\Triv$ if $q=-1$. It makes commutative the diagram
$$
\diagram{
\RepHq && \hidewidth\lmapr{12}{c}\hidewidth && \Grot\calH_n(q) \cr
\noalign{\smallskip}
& \mapse{e} && \mapne{d} & \cr
\noalign{\smallskip}
&& \Grot\calH_n(z)_Q && \cr
}
$$
where $d$ is the decomposition map and $c$ is the canonical map.
\endproclaim

\Proof.
Since $\calH_n(z)_Q$ is semisimple, the classes of Specht modules $S^\la_Q$, 
$\la\in\Par(n)$, for this algebra form a $\bbZ$-basis of $\Grot\calH_n(z)_Q$. 
Since these modules are absolutely irreducible, we have 
$\,\<S^\la_Q,S^\mu_Q\>=\de_{\la\mu}$ (the Kronecker symbol). In particular the 
bilinear form on $\Grot\calH_n(z)_Q$ is nondegenerate. The equality of 
dimensions in Lemma 2.1 can be restated by saying that
$$
\<\Nk,\Mk\>=\<N_Q,M_Q\>
$$
for any two lattices $M,N\in\Rep$ if $q\ne-1$ and for $M,N\in\Triv$ if $q=-1$. 

If $M,M'\in\Rep$ are such that $\Mk\cong\Mk'$ as $\calH_n(q)$-modules, 
and if moreover $M,M'\in\Triv$ when $q=-1$, then it follows from the displayed 
equality that $\<X,M_Q\>=\<X,M'_Q\>$ for each permutation $\calH_n(z)_Q$-module 
$X$, i.e. a module induced from the trivial 1-dimensional representation of a 
parabolic subalgebra of $\calH_n(z)_Q$. Since the classes of permutation modules 
form a $\bbZ$-basis of $\Grot\calH_n(z)_Q$, we conclude that $M_Q\cong M'_Q$ 
in this case.

This shows that the map $e$ is well-defined on the elements $[M_\bbk]$. By Lemma 
2.7 each indecomposable $\calH_n(q)$-module with a 1-dimensional source is 
isomorphic to $M_\bbk$ for a suitable choice of $M$. Hence $e$ is well-defined 
on the semigroup of positive elements in $\RepHq$. Since both $?\ot_O\bbk$ and 
$?\ot_OQ$ are additive functors, we have $e(a+b)=e(a)+e(b)$ for any pair of 
positive elements $a,b\in\RepHq$. It follows that $e$ extends to a group 
homomorphism on the whole $\RepHq$. Commutativity of the diagram is clear 
from the definition of $d$\/ in section 2.
\endproof

\setitemsize(iii)
\proclaim
Proposition 4.2.
There exist group homomorphisms\/ $\ch:\RepHq\to\Sym_n$ and 
$\psi:\Sym_n\to\Grot\calH_n(q)$ with the following properties\/{\rm:}

\item(i)
$\psi\circ\ch$ is equal to the canonical map $c:\RepHq\to\Grot\calH_n(q),$

\item(ii)
$\psi(h_\la e_\mu)=[\calH_n(q)\ot_{\calH_{\lamu}(q)}\bbk_\lamu]$ 
for each pair $(\lamu)\in\Par^2(n),$

\item(iii)
if $q\ne-1$ then $\ch([\calH_n(q)\ot_{\calH_{\lamu}(q)}\bbk_\lamu])
=h_\la e_\mu\,$ for each pair $(\lamu)\in\Par^2(n),$

\item(iv)
$\ch([\calH_n(q)\ot_{\calH_\la(q)}\bbk\triv])=h_\la$ for each $\la\in\Par(n),$

\item(v)
$\<\ch(a),\ch(b)\,\>=\<a,b\>$ for all $\,a,b\in\RepHq$.

\item(vi)
$\ch$\/ maps the subgroup $\,\TrivHq\,$ isomorphically onto $\,\Sym_n${\rm,}

\endproclaim

\Proof.
By the semisimple case recalled in section 1 there is a canonical isomorphism 
$$
\Grot\calH_n(z)_Q\cong\Sym_n.
$$
Composing the group homomorphisms $d$ and $e$ of Lemma 4.1 with the previous 
isomorphism, we obtain $\psi$ and $\ch$, respectively. The commutative 
diagram in Lemma 4.1 yields (i).

Now take $M=M^\lamu$, i.e. $M=\calH_n(z)\ot_{\calH_{\la,\mu}(z)}O_\lamu\in\Rep$ 
(see section 2). We have $d([M_Q])=[\Mk]$ by the definition of $d$. If 
$q\ne-1$ then $e([\Mk])=[M_Q]$ by the definition of $e$. If $\mu=0$, so that 
$M=M^\la\in\Triv$, then $e([\Mk])=[M_Q]$ for any $q$. Since 
$$
[M_Q]=[\calH_n(z)_Q\ot_{\calH_{\la,\mu}(z)_Q}Q_\lamu]\in\Grot\calH_n(z)_Q
$$ 
corresponds to $h_\la e_\mu\in\Sym_n$ by Lemma 1.1, we get (ii)--(iv).

By Lemma 2.7 the group $\RepHq$ is generated by the classes $[\Nk]$ for all 
lattices $N\in\Rep$, and when $q=-1$ it suffices to use only the lattices 
$N\in\Triv$. Since $e([\Nk])=[N_Q]$ for such lattices, it follows from 
Lemma 2.1 that
$$
\<e(a),e(b)\>=\<a,b\>
$$
for all $a,b\in\RepHq$. Since the isomorphism $\Sym_n\cong\Grot\calH_n(z)_Q$ 
is isometric with respect to the scalar products defined on these groups, we 
get (v).

The trivial source indecomposable $\calH_n(q)$-modules are precisely the Young 
modules $Y^\la$ parametrized by partitions $\la\in\Par(n)$. These modules were 
described by Dipper and James in \cite{Dip-J89, Lemma 2.5}. For each $\la$ 
the ``permutation" module
$$
\calH_n(q)\ot_{\calH_\la(q)}\bbk\triv
$$
has $Y^\la$ as its direct summand of multiplicity 1. All other direct summands 
of this permutation module are the Young modules $Y^\nu$ with $\nu>\la$ with 
respect to the dominance order on partitions.

The isomorphism classes of Young modules $Y^\la$ form a $\bbZ$-basis of 
$\TrivHq$. It follows that the classes of permutation modules also form a 
$\bbZ$-basis of $\TrivHq$. According to (iv) $\ch$ maps this basis of $\TrivHq$ 
to the $\bbZ$-basis $\{h_\la\mid\la\in\Par(n)\}$ of $\Sym_n$. This entails (vi).
\endproof

\Remark.
Actually the fact that the permutation $\calH_n(q)$-modules do not have any 
indecomposable direct summands other than the Young modules was not proved in 
\cite{Dip-J89} for arbitrary $q$. The isomorphism classes of those 
indecomposable summands are in a bijection with the isomorphism classes of 
simple modules for the $q$-Schur algebra $S_q(k,n)$ when $k\ge n$. In a later 
paper Dipper and James showed that the latter modules are parametrized by 
partitions of $n$ \cite{Dip-J91, Th. 8.8}. This settled the question about 
the direct summands of permutation modules (see Donkin \cite{Don, 4.4}).
\endremark

Let us return to consideration of the Hecke symmetry $R$. Put $\calH_n=\calH_n(q)$.

\proclaim
Lemma 4.3.
Let $f$ be the Hilbert series of the $R$-symmetric algebra $\bbS(V,R)$. Then
$$
\dim\Hom_{\calH_n}(V^{\ot n}\!,\calH_n\ot_{\calH_\nu}\bbk\triv)=\fhat(h_\nu)
=\smash{\sum_{\la\in\Par(n)}}\fhat(m_\la)\<h_\la,h_\nu\>
$$
for each $\nu\in\Par(n)$.
\endproclaim

\Proof.
By the Frobenius reciprocity
$$
\Hom_{\calH_n}(V^{\ot n}\!,\calH_n\ot_{\calH_\nu}\bbk\triv)\cong
\Hom_{\calH_\nu}(V^{\ot n}\!,\bbk\triv)\cong
\Homk(V^{\ot n}\!/\Si_\nu(V^{\ot n}),\bbk)
$$
where $\Si_\nu(V^{\ot n})=\sum_{i\in\calI_\nu}(T_i-q)V^{\ot n}$. If 
$\nu=(\nu_1,\ldots,\nu_k)$, then
$$
V^{\ot n}/\,\Si_\nu(V^{\ot n})\cong\bbS_{\nu_1}(V,R)\ot\ldots\ot\bbS_{\nu_k}(V,R)
$$
by Lemma 3.1. Since $\,\dim\bbS_{\nu_i}(V,R)=\fhat(h_{\nu_i})$, we get
$$
\dim V^{\ot n}\!/\Si_\nu(V^{\ot n})
=\prod_{i=1}^k\fhat(h_{\nu_i})=\fhat(h_{\nu_1}\cdots h_{\nu_k})=\fhat(h_\nu),
$$
which yields the first equality in the statement of Lemma 4.3. The second 
equality follows from the fact that the two bases $\{m_\la\mid\la\in\Par(n)\}$ 
and $\{h_\la\mid\la\in\Par(n)\}$ of the group $\Sym_n$ are dual to each other 
with respect to the scalar product on $\Sym_n$ \cite{Mac, Ch.~I, (4.5)}. Hence 
$u=\sum_{\la\in\Par(n)}\<h_\la,u\>m_\la$ for each $u\in\Sym_n$, and so this 
formula can be used for $u=h_\nu$.
\endproof

\proclaim
Lemma 4.4.
If\/ $f=\prod\limits_{j=1}^{r_1}(1+\be_jt)\cdot
\prod\limits_{i=1}^{r_0}(1-\al_it)^{-1},$ then
$$
\fhat(h_\nu)=\sum_{(\lamu)\in\Par^2(n)}\Nlamunu\,m_\la(\al)m_\mu(\be),
\qquad\nu\in\Par(n),
$$
where we put $u(\al)=u(\al_1,\ldots,\al_{r_0}),$ 
$u(\be)=u(\be_1,\ldots,\be_{r_1})$ for each $u\in\Sym$ and
$\Nlamunu$ is the number of all pairs of nonnegative integer matrices
$$
A=(a_{li})_{1\le l\le\ell(\nu),\,1\le i\le \ell(\la)},\qquad 
B=(b_{lj})_{1\le l\le\ell(\nu),\,1\le j\le \ell(\mu)}
$$
such that $B$ has only entries equal to $0$ or $1,$
$$
\displaylines{
{\textstyle\sum\limits_{l=1}^{\ell(\nu)}}\,a_{li}=\la_i
\ \ \hbox{\rm for each $i=1,\ldots,\ell(\la)$},\qquad
{\textstyle\sum\limits_{l=1}^{\ell(\nu)}}\,b_{lj}=\mu_j
\ \ \hbox{\rm for each $j=1,\ldots,\ell(\mu)$},\cr
{\rm and}\quad{\textstyle\sum\limits_{i=1}^{\ell(\la)}}\,a_{li}
+{\textstyle\sum\limits_{j=1}^{\ell(\mu)}}\,b_{lj}=\nu_l
\ \ \hbox{\rm for each $l=1,\ldots,\ell(\nu)$}.\cr
}
$$
Also{\rm,} $\,\Nlamunu=\<h_\la e_\mu,h_\nu\>$.
\endproclaim

\Proof.
Recall that $\,\prod(1-\al_it)^{-1}=\sum h_p(\al)t^p\,$ and
$\,\prod(1+\be_jt)=\sum e_p(\be)t^p$ by the definitions of the complete and 
elementary symmetric functions. For each integer $p\ge0$ the coefficient 
of $t^p$ in $f$ is
$$
\fhat(h_p)=\sum_{i=0}^ph_i(\al)e_{p-i}(\be),
$$
i.e. $\fhat(h_p)$ equals the sum of all monomials 
$\al_1^{a_1}\cdots\al_{r_0}^{a_{r_0}}\be_1^{b_1}\cdots\be_{r_1}^{b_{r_1}}$ 
with integer exponents $a_i\ge0$ and $b_j\in\{0,1\}$ such that 
$\sum a_i+\sum b_j=p$. Hence
$$
\fhat(h_\nu)=\prod_{l=1}^k\fhat(h_{\nu_l})=\sum_{A,B}\ \prod_{l=1}^k
\bigl(\al_1^{a_{l1}}\cdots\al_{r_0}^{a_{lr_0}}
\be_1^{b_{l1}}\cdots\be_{r_1}^{b_{lr_1}}\bigr)
$$
where $k=\ell(\nu)$ and the sum runs over all pairs of nonnegative 
integer matrices
$$
A=(a_{li})_{1\le l\le k,\,1\le i\le r_0},\qquad 
B=(b_{lj})_{1\le l\le k,\,1\le j\le r_1}
$$
such that $B$ has only entries equal to 0 or 1 and 
$$
{\textstyle\sum\limits_{i=1}^{r_0}}\,a_{li}
+{\textstyle\sum\limits_{j=1}^{r_1}}\,b_{lj}=\nu_l\quad
\hbox{for each $l=1,\ldots,k$}.
$$

For $\la,\mu\in\Par$ we have $m_\la(\al)=0$ whenever $\ell(\la)>r_0$ and 
$m_\mu(\be)=0$ whenever $\ell(\mu)>r_1$. If $\ell(\la)\le r_0$ and 
$\ell(\mu)\le r_1$, then each occurrence of the monomial 
$\al_1^{\la_1}\cdots\al_{r_0}^{\la_{r_0}}\be_1^{\mu_1}\cdots\be_{r_1}^{\mu_{r_1}}$ 
in the previous expression for $\fhat(h_\nu)$ corresponds to a pair of matrices 
$A,B$ used in the definition of $\Nlamunu$. Thus $\Nlamunu$ is the total 
number of such occurrences. Since $\fhat(h_\nu)$ is a symmetric function of 
$\al_1,\ldots,\al_{r_0}$ and a symmetric function of $\be_1,\ldots,\be_{r_1}$, 
it is a $\bbZ$-linear combination of the products of monomial symmetric 
functions in these sets of elements, and we obtain the desired formula. 

Next, writing out the product $h_\la e_\mu=h_{\la_1}\cdots h_{\la_{\ell(\la)}}
e_{\mu_1}\cdots e_{\mu_{\ell(\mu)}}$ as a sum of monomials in the 
indeterminates $x_1,x_2,\ldots\,$, we see that the coefficient of 
$x_1^{\nu_1}\cdots x_k^{\nu_k}$ in $h_\la e_\mu$ equals $\Nlamunu$ too. Hence
$$
h_\la e_\mu=\sum_{\rho\,\in\,\Par(n)}N_{(\lamu),\mskip1mu\rho}\,m_\rho\quad
\hbox{for $(\lamu)\in\Par^2(n)$},
$$
and it follows that $\,\<h_\la e_\mu,h_\nu\>=\Nlamunu\,$.
\endproof

We will write $\,\ch(X)=\ch([X])\,$ for each finite dimensional left 
$\calH_n$-module $X$ where $\,\ch\,$ is the map of Proposition 4.2. The property
$$
\<\ch(X),\ch(Y)\>=\<X,Y\>
$$
of this map enables us to determine $\ch(V^{\ot n})$:

\proclaim
Theorem 4.5.
Let $R:V^{\ot2}\to V^{\ot2}$ be a Hecke symmetry of birank $(r_0,r_1)$ with 
parameter $q$. Suppose that $R$ satisfies the $1$-dimensional source 
condition and
$$
H_{\bbS(V,R)}={\tprod_{j=1}^{r_1}}(1+\be_jt)\cdot
{\tprod_{i=1}^{r_0}}(1-\al_it)^{-1}.
$$
Then \ $\displaystyle\ch(V^{\ot n})
=\sum_{(\lamu)\in\Par^2(n)}\,m_\la(\al)m_\mu(\be)h_\la e_\mu\,$. In particular{\rm,}
$$
[V^{\ot n}]=\sum_{(\lamu)\in\Par^2(n)}\,m_\la(\al)m_\mu(\be)
[\calH_n(q)\ot_{\calH_{\lamu}(q)}\bbk_\lamu]
$$
in the Grothendieck group\/ $\Grot\calH_n(q)$.
\endproclaim

\Proof.
By Lemmas 4.3 and 4.4
$$
\<V^{\ot n},\calH_n\ot_{\calH_\nu}\bbk\triv\>=\fhat(h_\nu)
=\sum_{(\lamu)\in\Par^2(n)}m_\la(\al)m_\mu(\be)\<h_\la e_\mu,h_\nu\>
$$
for each $\nu\in\Par(n)$. Also, $\,\<\ch(V^{\ot n}),h_\nu\>=\<V^{\ot n},
\calH_n\ot_{\calH_\nu}\bbk\triv\>\,$ 
by (iv) and (v) of Proposition 4.2. This means that
$$
\ch(V^{\ot n})-\sum_{(\lamu)\in\Par^2(n)}\,m_\la(\al)m_\mu(\be)h_\la e_\mu
$$
is orthogonal to all functions $h_\nu$, $\nu\in\Par(n)$, which generate the 
whole group $\Sym_n$. Nondegeneracy of the scalar product entails the required 
formula for $\ch(V^{\ot n})$. The final equality is obtained then by 
applying the map $\psi:\Sym_n\to\Grot\calH_n$ described in Proposition 4.2.
\endproof

\proclaim
Corollary 4.6.
If $n$ is such that the Hecke algebra $\calH_n(q)$ is semisimple{\rm,} 
then{\rm,} as an $\calH_n(q)$-module{\rm,}
$$
V^{\ot n}\cong\bigoplus_{(\lamu)\in\Par^2(n)}
\bigl(\calH_n(q)\ot_{\calH_{\lamu}(q)}\bbk_\lamu\bigr)^{m_\la(\al)m_\mu(\be)}.
$$
\endproclaim

\Proof.
In this case the three maps $c,d,e$ of Lemma 4.1 are bijective. Hence so too 
is the map $\ch:\RepHq\to\Sym_n$. The two modules in the statement are 
isomorphic since they have the same image in $\Sym_n$.
\endproof

\section
5. The Hilbert series of the intertwining algebras

Let $K$ be any commutative ring. Denote by $\UK$ the multiplicative subgroup 
of the ring $K[[t]]$ consisting of all formal power series with constant term 
equal to 1. In other words, $\UK=1+\frm$ where $\frm$ is the ideal of $K[[t]]$ 
generated by $t$. There is a well known $\la$-ring structure on $\UK$ 
\cite{At-T69, Lemma 1.1}. Addition in this ring is given by the usual 
multiplication of power series, while the ring multiplication is another 
binary operation $\circ$ which has the property that $(1+at)\circ(1+bt)=1+abt$ 
for $a,b\in K$. However, Theorem 5.5 makes use of a different multiplication 
$\diamond$ which gives an isomorphic ring structure on $\UK$ and satisfies
$$
(1+at)\diamond(1+bt)=(1-at)^{-1}\diamond(1-bt)^{-1}=(1-abt)^{-1}.
$$

Now we will give a formal definition of this binary operation using the theory 
of symmetric functions. For each $n\ge0$ extend the scalar product on $\Sym_n$ 
to a symmetric $K$-bilinear form on the $K$-module $\Sym_{n,K}=K\ot_{\bbZ}\Sym_n$. 
Since this bilinear form induces a bijection 
$$
\Sym_{n,K}\mapr{}\Hom_K(\Sym_{n,K},K)\cong\HomZ(\Sym_n,K),
$$
to each $f\in\UK$ there corresponds a uniquely determined element 
$\xi_n(f)\in\Sym_{n,K}$ such that
$$
\fhat(u)=\<\xi_n(f),u\>\quad\hbox{for all $\,u\in\Sym_n$}
$$
where $\fhat:\Sym\to K$ is the ring homomorphism defined in section 3 and we 
identify elements $u\in\Sym_n$ with their images $1\ot u$ in $\Sym_{n,K}$. Put
$$ 
f\diamond g=\sum_{n=0}^\infty\,\<\xi_n(f),\xi_n(g)\>\,t^n\quad
\hbox{for $\,f,g\in\UK$}.
$$
Denote by $c_n(F)$ the coefficient of $t^n$ in a formal power series $F$. 
Considering the extensions of $\fhat$, $\ghat$ to $K$-algebra homomorphisms 
$K\ot_\bbZ\Sym\to K$, we have
$$
c_n(f\diamond g)=\<\xi_n(f),\xi_n(g)\>=\fhat\bigl(\xi_n(g)\bigr)
=\ghat\bigl(\xi_n(f)\bigr).
$$

\proclaim
Lemma 5.1.
The formula $\xi(f)=\sum_{n=0}^\infty\xi_n(f)t^n$ defines a group homomorphism 
$\xi:\UK\to U_{\!K\ot_\bbZ{}\Sym}$.
\endproclaim

\Proof.
Since $\fhat(1)=1$, it follows that $\xi_0(f)=1$, i.e. the map $\xi$ has values 
in the group $U_{\!K\ot_\bbZ{}\Sym}$. We have to show that $\xi(fg)=\xi(f)\xi(g)$ 
for all $f,g\in \UK$. Recall the Hopf algebra structure on $\Sym$ (see 
\cite{Mac, Ch. I, Example 5.25}). Let
$$
\De:\Sym\to\Sym\ot_\bbZ\Sym
$$
be the comultiplication in $\Sym$ and $m_K:K\ot_\bbZ K\to K$ the multiplication 
in $K$. The formula
$$
\psi=m_K\circ(\fhat\ot\ghat)\circ\De
$$
defines a ring homomorphism $\Sym\to K$. Since 
$\De(h_n)=\sum_{i=0}^nh_i\ot h_{n-i}$, we have
$$
\psi(h_n)=\sum_{i=0}^n\fhat(h_i)\ghat(h_{n-i})
=\sum_{i=0}^nc_i(f)c_{n-i}(g)=c_n(fg)
$$
for each $n\ge0$. It follows that $\psi=\phat$ where $p=fg$. Hence
$$
\phat(u)
=\bigl(m_K\circ(\fhat\ot\ghat)\circ\De\bigr)(u)
=\sum_{i=0}^n\<\xi_i(f)\ot\xi_{n-i}(g),\De(u)\>
=\sum_{i=0}^n\<\xi_i(f)\xi_{n-i}(g),u\>
$$
for each $u\in\Sym_n$. The last equality in this formula follows from the fact 
that $\De$ is the adjoint of the multiplication in $\Sym$ with respect to the 
scalar product. The middle equality is true because 
$\De(u)\in\sum_{j=0}^n\Sym_j\ot_\bbZ\Sym_{n-j}$ and
$$
\<\xi_i(f)\ot\xi_{n-i}(g),u'\ot u''\>
=\<\xi_i(f),u'\>\<\xi_{n-i}(g),u''\>=\cases{\fhat(u')\ghat(u'') & when $j=i$\cr
\noalign{\smallskip}
0 & when $j\ne i$}
$$
for any pair of elements $u'\in\Sym_j$ and $u''\in\Sym_{n-j}$.

We conclude that $\xi_n(p)=\sum_{i=0}^n\xi_i(f)\xi_{n-i}(g)$ for each 
$n\ge0$, which amounts to the desired equality of formal power series 
$\xi(p)=\xi(f)\xi(g)$.
\endproof

\proclaim
Lemma 5.2.
Let $f,g_1,g_2\in\UK$. Then $f\diamond(g_1g_2)=(f\diamond g_1)(f\diamond g_2)$.
\endproclaim

\Proof.
Since $\xi_n(g_1g_2)=\sum_{i=0}^n\xi_i(g_1)\xi_{n-i}(g_2)$ by Lemma 5.1, we 
have
$$
\fhat\bigl(\xi_n(g_1g_2)\bigr)
=\sum_{i=0}^n\fhat\bigl(\xi_i(g_1)\bigr)\fhat\bigl(\xi_{n-i}(g_2)\bigr)\quad
\hbox{for each $n\ge0$}.
$$
Thus the two formal power series in question have equal coefficients.
\endproof

\proclaim
Lemma 5.3.
Let $a\in K$. If $g=(1-at)^{-1},$ then $f\diamond g=f(at)$. If $g=1+at,$ then 
$f\diamond g=f(-at)^{-1}$.
\endproclaim

\Proof.
It follows from Lemma 5.2 that $f\diamond g^{-1}=(f\diamond g)^{-1}$ for all 
$g\in\UK$. So the second assertion of Lemma 5.3 is equivalent to the first 
one. Let $g=(1-at)^{-1}$. Then $\ghat(h_n)=a^n$ for all $n\ge0$. Hence 
$\ghat(h_\nu)=a^{|\nu|}$ for all $\nu\in\Par$. Since the two $\bbZ$-bases 
$\{h_\nu\}$ and $\{m_\nu\}$ of $\Sym$ are dual to each other, we get
$$
\xi_n(g)=a^n\sum_{\la\in\Par(n)}m_\la=a^nh_n.
$$
Therefore, $c_n(f\diamond g)=\fhat\bigl(\xi_n(g)\bigr)=a^n\fhat(h_n)=a^nc_n(f)$, 
which is the coefficient of $t^n$ in $f(at)$.
\endproof

\Remark.
If $p=f\diamond g$, then $\phat=m_K\circ(\fhat\ot\ghat)\circ\De^*$ where
$\De^*:\Sym\to\Sym\ot_\bbZ\Sym$ is the ring homomorphism discussed in 
\cite{Mac, Ch. I, Example 7.20}. Associativity of $\diamond$ follows from 
coassociativity of $\De^*$.
\endremark

\proclaim
Lemma 5.4.
Let $f$ be the Hilbert series of the $R$-symmetric algebra\/ $\bbS(V,R)$. If 
$R$ satisfies the $1$-dimensional source condition, then 
$\xi_n(f)=\ch(V^{\ot n})$ where\/ $\ch$ is the map of Proposition\/ {\rm4.2}.
\endproclaim

\Proof.
Here $K=\bbZ$, and so $\xi_n(f)\in\Sym_n$. As we have seen in the proof of 
Theorem 4.5,
$$
\<\ch(V^{\ot n}),h_\nu\>=\<V^{\ot n},\calH_n\ot_{\calH_\nu}\bbk\triv\>
=\fhat(h_\nu)
$$
for all $\nu\in\Par(n)$. Hence $\<\ch(V^{\ot n}),u\>=\fhat(u)$ for all 
$u\in\Sym_n$, and the conclusion follows from the definition of $\xi_n(f)$.
\endproof

Suppose now that $V'$ is a second finite dimensional vector space over the same 
field $\bbk$, and $R'$ is a Hecke symmetry on $V'$ with the same parameter $q$ 
as the Hecke symmetry $R$ on $V$. We will define the algebra $A(R',R)$ (cf. 
\cite{Hai02} and \cite{Sk19}).

Consider the tensor algebra of the vector space ${\Homk(V,V')}^{\mskip-1mu*}$ 
dual to $\Homk(V,V')$. Its homogeneous component of degree $n$ admits the 
following realization via canonical $\bbk$-linear bijections
$$
\bbT_n\bigl({\Homk(V,V')}^{\mskip-2mu*}\bigr)\cong
\bigl(\Homk(V,V')^{\ot n}\bigr)^{\mskip-2mu*}\cong
\Homk\bigl(V^{\ot n},{V'}^{\mskip2mu\ot n}\bigr)^{\mskip-2mu*}.
$$

Denote by $\calR$ the linear operator on the vector space 
$\Homk\bigl(V^{\ot2},{V'}^{\mskip2mu\ot2}\bigr)^{\mskip-2mu*}$ such that the 
dual operator $\calR^*$ on $\Homk(V^{\ot2},{V'}^{\mskip2mu\ot2})$ is 
defined by the formula
$$
\calR^*(h)={(R')}^{-1}\circ h\circ R,\qquad 
h\in\Homk\bigl(V^{\ot2},{V'}^{\mskip2mu\ot2}\bigr).
$$
In other words, $\calR={({R'}^{\mskip2mu*})}^{\mskip-2mu-1}\ot R$ under 
the canonical identification 
$$
\Homk\bigl(V^{\ot2},{V'}^{\mskip2mu\ot2}\bigr)^{\mskip-2mu*}
\cong\bigl({V'}^{\mskip2mu\ot2}\bigr)^{\mskip-2mu*}\ot V^{\ot2}.
$$

The algebra $A(R',R)$ is defined as the factor algebra of 
$\bbT\bigl({\Homk(V,V')}^{\mskip-1mu*}\bigr)$ by the homogeneous 
ideal $I$ generated by
$$
I_2=\Im\,(\calR-\Id)\sbs
\Homk\bigl(V^{\ot2},{V'}^{\mskip2mu\ot2}\bigr)^{\mskip-2mu*}\cong
\bbT_2\bigl({\Homk(V,V')}^{\mskip-2mu*}\bigr).
$$
Under the canonical pairings between 
$\bbT_n\bigl({\Homk(V,V')}^{\mskip-2mu*}\bigr)$ and 
$\Homk\bigl(V^{\ot n},{V'}^{\mskip2mu\ot n}\bigr)$ we have
$$
\eqalign{
I_2^\perp&{}=\{h\in\Homk\bigl(V^{\ot2},{V'}^{\mskip2mu\ot2}\bigr)\mid\calR^*(h)=h\}\cr
&{}=\{h\in\Homk\bigl(V^{\ot2},{V'}^{\mskip2mu\ot2}\bigr)\mid h\circ R=R'\circ h\}.
}
$$
for $n=2$. In degree $n>2$ the homogeneous component of the ideal $I$ is
$$
I_n=\sum_{i=1}^{n-1}\,\bbT_{i-1}\bigl({\Homk(V,V')}^{\mskip-2mu*}\bigr)
\ot I_2\ot\bbT_{n-i-1}\bigl({\Homk(V,V')}^{\mskip-2mu*}\bigr).
$$
Consider $V^{\ot n}$ and ${V'}^{\mskip2mu\ot n}$ as $\calH_n(q)$-modules with 
respect to representations arising from the Hecke symmetries $R$ and $R'$, 
respectively. It follows then that
$$
\eqalign{
I_n^\perp=\{h\in\Homk\bigl(V^{\ot n},{V'}^{\mskip2mu\ot n}\bigr)\mid
hT_i=T_ih\,\hbox{ for all }\,i=1,\ldots n-1\}\qquad\qquad&\cr
=\Hom_{\calH_n}\!\bigl(V^{\ot n},{V'}^{\mskip2mu\ot n}\bigr).&
}
$$
This equality holds for all $n\ge0$ since $I_n=0$ when $n\le1$. Hence
$$
A_n(R',R)=\bbT_n\bigl({\Homk(V,V')}^{\mskip-2mu*}\bigr)/I_n\cong
\Hom_{\calH_n}\!\bigl(V^{\ot n},{V'}^{\mskip2mu\ot n}\bigr)^{\mskip-2mu*}\quad
\hbox{for all $n$}.
$$

Note that $A(R',R)=A(R)$ when $R'=R$. Another graded algebra $E(R',R)$ related 
to $A(R',R)$ is defined as the factor algebra of 
$\bbT\bigl({\Homk(V,V')}^{\mskip-1mu*}\bigr)$ by the ideal generated by 
$\Ker\,(\calR-\Id)$.

\proclaim
Theorem 5.5.
Let $R,$ $R'$ be Hecke symmetries on finite dimensional vector spaces $V,$ 
$V'$ with the same parameter $q$ of the Hecke relation. Suppose that both $R$ 
and $R'$ satisfy the $1$-dimensional source condition. Let
$$ 
f={\tprod_{j=1}^{r_1}}(1+\be_jt)\cdot
{\tprod_{i=1}^{r_0}}(1-\al_it)^{-1},\qquad
g={\tprod_{j=1}^{r'_1}}(1+\be'_jt)\cdot
{\tprod_{i=1}^{r'_0}}(1-\al'_it)^{-1}
$$
be the Hilbert series of the algebras $\bbS(V,R)$ and $\bbS(V'\!,R')$. Then 
the Hilbert series of the algebra $A(R',R)$ equals
$$
{\tprod_{i=1}^{r_1}}{\tprod_{j=1}^{r'_0}}(1+\be_i\al'_jt)\cdot
{\tprod_{i=1}^{r_0}}{\tprod_{j=1}^{r'_1}}(1+\al_i\be'_jt)\cdot
{\tprod_{i=1}^{r_0}}{\tprod_{j=1}^{r'_0}}(1-\al_i\al'_jt)^{-1}\cdot
{\tprod_{i=1}^{r_1}}{\tprod_{j=1}^{r'_1}}(1-\be_i\be'_jt)^{-1}.
$$
\endproclaim

\Proof.
The isomorphism classes of $\calH_n(q)$-modules $V^{\ot n}$ and 
${V'}^{\mskip2mu\ot n}$ define elements of $\RepHq$. We have $\xi_n(f)=\ch(V^{\ot n})$ 
and $\xi_n(g)=\ch({V'}^{\mskip2mu\ot n})$ by Lemma 5.4. Now
$$
\dim A_n(R',R)=\dim\Hom_{\calH_n(q)}\bigl(V^{\ot n},{V'}^{\mskip2mu\ot n}\bigr)
=\<V^{\ot n},{V'}^{\mskip2mu\ot n}\>.
$$
Proposition 4.2(v) and Lemma 5.4 yield
$$
\dim A_n(R',R)=\<\ch(V^{\ot n}),\,\ch({V'}^{\mskip2mu\ot n})\>=\<\xi_n(f),\xi_n(g)\>.
$$
Hence the Hilbert series $\sum_{n\ge0}\bigl(\dim A_n(R',R)\bigr)t^n$ 
of the algebra $A(R',R)$ coincides with $f\diamond g$. Taking $K=\bbR$ and
applying Lemmas 5.2 and 5.3, we get
$$
f\diamond g
={\tprod_{i=1}^{r'_0}}f(\al'_it)\cdot{\tprod_{j=1}^{r'_1}}f(-\be'_jt)^{-1},
$$ 
which gives the required formula.
\endproof

\section
6. Even symmetries satisfying the trivial source condition

If $R$ is a Hecke symmetry of birank $(r_0,r_1)$ with $r_1=0$, then $R$ is said 
to have rank $r=r_0$. In the situation of Theorem 3.8 this happens precisely 
when the Hilbert series of $\La(V,R)$ is a polynomial of degree $r$, i.e. the 
algebra $\La(V,R)$ is finite dimensional with the grading of length $r$. In 
\cite{Gur90} Gurevich calls a closed Hecke symmetry of rank $r$ \emph{even}. 
We do not need the assumption of closedness, however. Put
$$
\Par(n,r)=\{\la\in\Par(n)\mid\ell(\la)\le r\}.
$$

\setitemsize(iii)
\proclaim
Theorem 6.1.
Suppose that $R$ is a Hecke symmetry of rank $r$ satisfying the trivial source 
condition. Let\/ $\tprod_{i=1}^r(1+\al_it)$ be the Hilbert series of 
$\La(V,R)$. For each $n\ge0$

\item(i)
$V^{\ot n}\cong\bigoplus_{\la\in\Par(n)}
\bigl(\calH_n(q)\ot_{\calH_\la(q)}\bbk\triv\bigr)^{m_\la(\al)}$ as 
$\calH_n(q)$-modules{\rm,}

\smallskip
\item(ii)
the algebra $A_n(R)^*$ is Morita equivalent to the $q$-Schur algebra $S_q(r,n),$

\smallskip
\item(iii)
$\{[V^\la]\mid\la\in\Par(n,r)\}$ is a $\bbZ$-basis of\/ $\Grot_n(R)$.

\smallskip
There is an isomorphism of graded rings\/ $\Grot(R)\cong\Sym/I_r$ where $I_r$ 
is the ideal of\/ $\Sym$ with a $\bbZ$-basis\/ 
$\{s_\la\mid\la\in\Par,\ \ell(\la)>r\}$.
\endproclaim

\Proof.
Under present assumptions the $\calH_n(q)$-module $V^{\ot n}$ defines an 
element of the group $\TrivHq$, and so too does the $\calH_n(q)$-module in the 
right hand side of (i). By Theorem 4.5 
$$
\ch(V^{\ot n})=\sum_{\la\in\Par(n)}m_\la(\al)h_\la.
$$
Parts (iv) and (vi) of Proposition 4.2 show that the second module has the 
same image under the isomorphism of groups $\,\ch:\TrivHq\to\Sym_n$, and 
therefore the two modules are isomorphic. This is conclusion (i).

Recall that $A_n(R)^*\cong E_n=\End_{\calH_n(q)}V^{\ot n}$, while $S_q(r,n)$ 
is the endomorphism ring of the $\calH_n(q)$-module
$$
M(r,n)=\bigoplus_{\mu=(\mu_1,\ldots,\mu_r)}
\bigl(\calH_n(q)\ot_{\calH_\mu(q)}\bbk\triv\bigr)
$$
where the sum is taken over all weak compositions of $n$ with $r$ components. 
A \emph{weak composition} is allowed to have zero components $\mu_i$, and the 
corresponding parabolic subalgebra $\calH_\mu(q)$ is the same as one would get 
by removing from $\mu$ all its zero components. If $\mu''\in\Par(n,r)$ is the 
partition obtained from $\mu$ by rearranging its nonzero components in 
decreasing order, then
$$ 
\calH_n(q)\ot_{\calH_\mu(q)}\bbk\triv\cong
\calH_n(q)\ot_{\calH_{\mu''}(q)}\bbk\triv\,.
$$
On the other hand, $m_\la(\al)\ge0$ for each $\la\in\Par(n)$, and $m_\la(\al)=0$ 
if and only if $\ell(\la)>r$. Therefore nonzero summands in the right hand side 
of (i) correspond precisely to partitions $\la\in\Par(n,r)$. We see that the 
$\calH_n(q)$-modules $V^{\ot n}$ and $M(r,n)$ have the same set of isomorphism 
classes of indecomposable direct summands. But then the endomorphism rings $E_n$ 
and $S_q(r,n)$ of the two modules are Morita equivalent. This entails (ii).

Although the direct sum decomposition of (i) in general does not come from an 
action of a torus, we can still use it to describe the structure of 
$A(R)$-comodules by means of a kind of weight spaces. Recall that right 
$A_n(R)$-comodules may be identified with left $E_n$-modules. Fix an 
isomorphism in (i). It gives a collection of $\calH_n(q)$-submodules
$$ 
M_\lai\sbs V^{\ot n}\qquad\hbox{with $\,\la\in\Par(n,r)$, 
$\,1\le i\le m_\la(\al)$}
$$
such that $V^{\ot n}=\bigoplus M_\lai$ and 
$M_\lai\cong\calH_n(q)\ot_{\calH_\la(q)}\bbk\triv$ for each pair $(\la,i)$. 
Denote by $\xi_\lai$ the projection onto $M_\lai$ with respect to this 
decomposition. Then
$$
\{\xi_\lai\mid\la\in\Par(n,r),\ 1\le i\le m_\la(\al)\}
$$
is a set of pairwise orthogonal idempotents in the ring $E_n$ with 
$\sum\xi_\lai=1$. If $X$ is any left $E_n$-module, then $X=\bigoplus X_\lai$ 
where $X_\lai=\xi_\lai X$. Put $X_\la=X_{\la,1}$ and define the 
\emph{character}\/ of $X$ by the formula
$$
\chR(X)=\sum_{\la\in\Par(n,r)}(\dim X_\la)\,m_\la(x_1,\ldots,x_r)
\in\Sym_n(r)\sbs\bbZ[x_1,\ldots,x_r]
$$
where $\Sym_n(r)$ is the group of all symmetric homogeneous polynomials of 
degree $n$ in $r$ indeterminates $x_1,\ldots,x_r$.

Since $M_\lai\cong M_{\la,1}$, the idempotent $\xi_\lai$ is conjugate to 
$\xi_{\la,1}$ by an inner automorphism of $E_n$. Hence $\dim X_\lai=\dim X_\la$ 
for all $i$ such that $1\le i\le m_\la(\al)$. In particular, it follows that 
$X=0$ whenever $\chR(X)=0$.

If $\,\,0\to X'\to X\to X''\to0\,\,$ is an exact sequence of finite dimensional 
left $E_n$-modules, then $\,\,\dim X_\la=\dim X'_\la+\dim X''_\la\,\,$ for each 
$\la$, whence 
$$
\chR(X)=\chR(X')+\chR(X'').
$$
So $\chR$ gives rise to a group homomorphism
$$
\Grot_n(R)\cong\Grot E_n\mapr{}\Sym_n(r)
$$
which will be denoted by the same symbol $\chR$.

It is well known that the Schur polynomials $s_\la(x_1,\ldots,x_r)$ with 
$\la\in\Par(n,r)$ form a $\bbZ$-basis of $\Sym_n(r)$. On the other hand, 
isomorphism classes of simple $S_q(n,r)$-modules are also parametrized by the 
set $\Par(n,r)$. The Morita equivalent algebra $E_n$ has the same number of 
simple modules equal to the cardinality of $\Par(n,r)$. In other words, 
$\Grot_n(R)$ and $\Sym_n(r)$ are free abelian groups of equal ranks. As will 
be proved separately in Lemma 6.2, $s_\la(x_1,\ldots,x_r)$ is the image of 
$V^\la$ under $\chR\mskip-2mu$. Hence $\chR$ is surjective, but then $\chR$ 
has to map the group $\Grot_n(R)$ isomorphically onto $\Sym_n(r)$, and (iii) 
is also clear.

The ring homomorphism $\ph:\Sym\to\Grot(R)$ defined in Proposition 3.2 sends 
$s_\la$ to $[V^\la]$ for each $\la\in\Par$. It follows from (iii) that $\ph$ 
maps the subgroup of $\Sym$ generated by $\{s_\la\mid\ell(\la)\le r\}$ 
isomorphically onto $\Grot(R)$. In particular, $\ph$ is surjective. If 
$\ell(\la)>r$, then $V^\la=0$ by Corollary 3.10. This can be also seen from 
the fact that $\,\chR(V^\la)=0\,$: assuming that $|\la|=n$, we have
$$
\chR(V^\la)=s_\la(x_1,\ldots,x_r)
=\sum_{\mu\in\Par(n)}K_{\la\mu}\mskip2mu m_\mu(x_1,\ldots,x_r).
$$
Recall that $K_{\la\mu}=0$ unless $\mu\le\la$ with respect to the dominance 
order. If $\mu\le\la$, then $\ell(\mu)\ge\ell(\la)>r$, in which case 
$m_\mu(x_1,\ldots,x_r)=0$. Thus all summands in the above expression vanish. 

It follows that $\,\Ker\ph=I_r\,$, and we are done.
\endproof

\proclaim
Lemma 6.2.
For each $\la\in\Par(n)$ we have $\chR(V^\la)=s_\la(x_1,\ldots,x_r)$.
\endproclaim

\Proof.
First we evaluate $\chR(\mskip2mu\Smu)$ for $\mu\in\Par(n)$. Recall that 
$\,\Smu\cong\bbk\triv\ot_{\calH_\mu(q)}V^{\ot n}$. If $\la\in\Par(n,r)$, then
$$
\Smula\cong\bbk\triv\ot_{\calH_\mu(q)}M_{\la,1}
\cong\bbk\triv\ot_{\calH_\mu(q)}\calH_n(q)\ot_{\calH_\la(q)}\bbk\triv
$$
since the $E_n$-module structure on $\Smu$ comes from the action of $E_n$ 
on $V^{\ot n}$. Hence $\dim\Smula$ is equal to the number $N_{\mu\la}$ 
of all $\frS_\mu\,$-$\,\frS_\la$ double cosets in the group $\frS_n$. There is 
also a combinatorial description of $N_{\mu\la}=N_{\la\mu}$ as the number of 
all nonnegative integer matrices of size $\ell(\la)\times\ell(\mu)$ having row 
sums $\la_i$, $\,i=1,\ldots,\ell(\la)$, and column sums $\mu_j$, 
$\,j=1,\ldots,\ell(\mu)$. By \cite{Mac, Ch. I, (6.7)}
$\,h_\mu= \sum_{\la\in\Par(n)}N_{\mu\la}\,m_\la$. If $\ell(\la)>r$, then 
$m_\la(x_1,\ldots,x_r)=0$. Therefore
$$
\chR(\,\Smu)=\sum_{\la\in\Par(n,r)}N_{\mu\la}\,m_\la(x_1,\ldots,x_r)
=h_\mu(x_1,\ldots,x_r).
$$
Let $\ph:\Sym\to\Grot(R)$ be the ring homomorphism defined in Proposition 3.2. 
Since $\ph(h_\mu)=[\,\Smu]$, the previous equality can be rewritten as
$$
\chR\bigl(\ph(u)\bigr)=u(x_1,\ldots,x_r)
$$
for $u=h_\mu$. Since the set $\{h_\mu\mid\mu\in\Par(n)\}$ generates the whole 
group $\Sym_n$, the formula above holds then for all $u\in\Sym_n$. Taking 
$u=s_\la$, we get the required conclusion.
\endproof

\Remark.
For each $\la\in\Par(n,r)$ it follows from Lemma 6.2 that $\dim\,(V^\la)_\la=1$ 
and $(V^\la)_\mu=0$ unless $\mu\le\la$. In this sense $\la$ is the highest 
weight of $V^\la$. Since the function $X\mapsto\dim X_\la$ is additive on exact 
sequences of $E_n$-modules, there is exactly one composition factor $L^\la$ of 
$V^\la$ with nonzero $\la$-weight space. Clearly $L^\la$ is also a module of 
highest weight $\la$. In particular, $L^\la\not\cong L^\mu_{\vphantom{\la}}$ 
whenever $\la,\mu\in\Par(n,r)$ and $\la\ne\mu$. This implies that 
$\{L^\la\mid\la\in\Par(n,r)\}$ is the full set of pairwise nonisomorphic 
simple $E_n$-modules.
\endremark

Let $H(R)$ stand for the Hopf envelope of $A(R)$. If $R$ and $R'$ are two 
closed Hecke symmetries with the same parameter $q$ and of the same birank 
$(r_0,r_1)$, then in the semisimple case, by a theorem of \PHHai \cite{Hai05, 
Th.  4.3}, there is a braided monoidal equivalence between the categories of 
right $H(R)$-comodules and right $H(R')$-comodules.

We are interested in the bialgebra version of this result. Let $R$ and $R'$ be 
not necessarily closed Hecke symmetries on vector spaces $V$ and $V'$, 
respectively, with the same parameter $q$ of the Hecke relation. Consider the 
$\calH_n(q)$-module structures on $V^{\ot n}$ and ${V'}^{\mskip2mu\ot n}$ 
arising from $R$ and $R'$. It was proved in \cite{Sk19, Th. 7.2} that there 
is a braided monoidal equivalence between the categories of right 
$A(R)$-comodules and right $A(R')$-comodules provided that for each $n>0$ the 
indecomposable $\calH_n(q)$-modules isomorphic to direct summands of 
$V^{\ot n}$ are the same as those isomorphic to direct summands of 
${V'}^{\mskip2mu\ot n}$.

If $R$ and $R'$ have the same birank $(r_0,r_1)$, then in the semisimple case 
the required condition on direct summands is satisfied since in each of the 
two $\calH_n(q)$-modules the simple submodules are precisely the Specht 
modules $S^\la$ with $\la\in\Ga(r_0,r_1)\cap\Par(n)$. Unfortunately we are not 
able to extend this result to the nonsemisimple case because Theorem 4.5 does 
not provide enough information on the module structure of the tensor powers. 
Therefore we have to restrict the class of Hecke symmetries in the following 
statement:

\proclaim
Theorem 6.3.
Let $R$ and $R'$ be Hecke symmetries with the same parameter $q,$ of the same 
rank $r${\rm,} and both satisfying the trivial source condition. Then there is 
a braided monoidal equivalence between the categories of right {\rm(}or 
left{\rm)} $A(R)$-comodules and $A(R')$-comodules.
\endproclaim

\Proof.
By Theorem 6.1 the $\calH_n(q)$-modules $V^{\ot n}$ and ${V'}^{\mskip2mu\ot n}$ 
have the same set of isomorphism classes of indecomposable direct summands. So 
\cite{Sk19, Th. 7.2} does apply.
\endproof

In conclusion we improve yet another result from \cite{Sk19}. Definitions of 
the algebras $A(R',R)$ and $E(R',R)$ have been recalled in section 5.

\proclaim
Theorem 6.4.
Let $R$ and $R'$ be Hecke symmetries with the same parameter $q,$ both 
satisfying the trivial source condition. Suppose that the algebras $\La(V,R)$ 
and $\La(V'\!,R')$ are Frobenius with the gradings of length $r$ and $r',$ 
respectively. Then the algebra $E(R',R)$ is Frobenius with the grading of 
length $rr',$ while $A(R',R)$ is Gorenstein of global dimension $rr'$.  
\endproclaim

\Proof.
The Hilbert series of the algebras $\La(V,R)$ and $\La(V'\!,R')$ are 
polynomials of degree $r$ and $r'$, respectively. We can write them as
$$
\tprod_{i=1}^r(1+\al_it)\quad{\rm and}\quad\tprod_{j=1}^{r'}(1+\al'_jt).
$$
The Hilbert series of $\bbS(V,R)$ and $\bbS(V'\!,R')$ are 
$\prod\mskip1mu(1-\al_it)^{-1}$ and $\prod\mskip1mu(1-\al'_jt)^{-1}$ by 
Theorem 3.8. For the graded algebras $A(R',R)$ and $E(R',R)$ it holds then
$$ 
H_{A(R'\!,R)}=\tprod_{i=1}^r\tprod_{j=1}^{r'}(1-\al_i\al'_jt)^{-1},\qquad
H_{E(R'\!,R)}=\tprod_{i=1}^r\tprod_{j=1}^{r'}(1+\al_i\al'_jt).
$$
The first formula here follows from Theorem 5.5. The second formula is a 
consequence of the relation $\,H_{E(R'\!,R)}(t)\,H_{A(R'\!,R)}(-t)=1\,$ proved 
in \cite{Sk19, Th. 6.2}.

Since $\dim\La_r(V,R)=\dim\La_{r'}(V'\!,R')=1$, we have 
$\prod\al_i=\prod\al'_j=1$. Thus $H_{E(R'\!,R)}$ is a polynomial of degree 
$rr'$ with the leading coefficient equal to 1. The conclusion follows then 
from \cite{Sk19, Th. 6.6}.
\endproof

\references
\nextref Ais-SW52
\auth{M.,Aissen;I.J.,Schoenberg;A.M.,Whitney}
\paper{On the generating functions of totally positive sequences. I}
\journal{J.~Analyse Math.}
\Vol{2}
\Year{1952}
\Pages{93-103}

\nextref At-T69
\auth{M.F.,Atiyah;D.O.,Tall}
\paper{Group representations, $\lambda$-rings and the $J$-homomorphism}
\journal{Topology}
\Vol{8}
\Year{1969}
\Pages{253-297}

\nextref Ber-R87
\auth{A.,Berele;A.,Regev}
\paper{Hook Young diagrams with applications to combinatorics and to representations of Lie superalgebras}
\journal{Adv. Math.}
\Vol{64}
\Year{1987}
\Pages{118-175}

\nextref Br-K03
\auth{J.,Brundan;J.,Kujawa}
\paper{A new proof of the Mullineux conjecture}
\journal{J.~Algebr. Combin.}
\Vol{18}
\Year{2003}
\Pages{13-39}

\nextref Dav00
\auth{A.A.,Davydov}
\paper{Totally positive sequences and $R$-matrix quadratic algebras}
\journal{J.~Math. Sci. (New York)}
\Vol{100}
\Year{2000}
\Pages{1871-1876}

\nextref Dip-J86
\auth{R.,Dipper;G.,James}
\paper{Representations of Hecke algebras of general linear groups}
\journal{Proc. London Math. Soc.}
\Vol{52}
\Year{1986}
\Pages{20-52}

\nextref Dip-J89
\auth{R.,Dipper;G.,James}
\paper{The q-Schur algebra}
\journal{Proc. London Math. Soc.}
\Vol{59}
\Year{1989}
\Pages{23-50}

\nextref Dip-J91
\auth{R.,Dipper;G.,James}
\paper{$q$-Tensor space and $q$-Weyl modules}
\journal{Trans. Amer. Math. Soc.}
\Vol{327}
\Year{1991}
\Pages{251-282}

\nextref Don
\auth{S.,Donkin}
\book{The $q$-Schur Algebra}
\publisher{Cambridge Univ. Press}
\Year{1998}

\nextref Don01
\auth{S.,Donkin}
\paper{Symmetric and exterior powers, linear source modules and representations of Schur superalgebras}
\journal{Proc. London Math. Soc.}
\Vol{83}
\Year{2001}
\Pages{647-680}

\nextref Du-PW91
\auth{J.,Du;B.,Parshall;J.,Wang}
\paper{Two-parameter quantum linear groups and the hyperbolic invariance of $q$-Schur algebras}
\journal{J.~London Math. Soc.}
\Vol{44}
\Year{1991}
\Pages{420-436}

\nextref Edr52
\auth{A.,Edrei}
\paper{On the generating functions of totally positive sequences. II}
\journal{J.~Analyse Math.}
\Vol{2}
\Year{1952}
\Pages{104-109}

\nextref Geck-P
\auth{M.,Geck;G.,Pfeiffer}
\book{Characters of Finite Coxeter Groups and Iwahori-Hecke Algebras}
\publisher{Clarendon Press}
\Year{2000}

\nextref Gia-LDW17
\auth{E.,Giannelli;K.J.,Lim;W.,O'Donovan;M.,Wildon}
\paper{On signed Young permutation modules and signed $p$-Kostka numbers}
\journal{J.~Group Theory}
\Vol{20}
\Year{2017}
\Pages{637-679}

\nextref Gur90
\auth{D.I.,Gurevich}
\paper{Algebraic aspects of the quantum Yang-Baxter equation\inRus}
\journal{Algebra i Analiz}
\Vol{2:4}
\Year{1990}
\Pages{119-148}
\etransl{Leningrad Math. J.}
\Vol{2}
\Year{1991}
\Pages{801-828}

\nextref Hai99
\auth{P.H.,Hai}
\paper{Poincar\'e series of quantum spaces associated to Hecke operators}
\journal{Acta Math. Vietnam}
\Vol{24}
\Year{1999}
\Pages{235-246}

\nextref Hai02
\auth{P.H.,Hai}
\paper{Realizations of quantum hom-spaces, invariant theory, and quantum determinantal ideals}
\journal{J.~Algebra}
\Vol{248}
\Year{2002}
\Pages{50-84}

\nextref Hai05
\auth{P.H.,Hai}
\paper{On the representation categories of matrix quantum groups of type $A$}
\journal{Vietnam J. Math.}
\Vol{33}
\Year{2005}
\Pages{357-367}

\nextref Hen
\auth{P.,Henrici}
\book{Applied and Computational Complex Analysis, Vol. $1$}
\publisher{Wiley}
\Year{1974}

\nextref Lar-T91
\auth{R.G.,Larson;J.,Towber}
\paper{Two dual classes of bialgebras related to the concepts of ``quantum group'' and ``quantum Lie algebra''}
\journal{Comm. Algebra}
\Vol{19}
\Year{1991}
\Pages{3295-3345}

\nextref Mac
\auth{I.G.,Macdonald}
\book{Symmetric Functions and Hall Polynomials}
second edition,
\publisher{Oxford Univ. Press}
\Year{1995}

\nextref Resh-TF89
\auth{N.Yu.,Reshetikhin;L.A.,Takhtajan;L.D.,Faddeev}
\paper{Quantization of Lie groups and Lie algebras\inRus}
\journal{Algebra i Analiz}
\Vol{1:1}
\Year{1989}
\Pages{178-206}
\etransl{Leningrad Math. J.}
\Vol{1}
\Year{1990}
\Pages{193-225}

\nextref Ser84
\auth{A.N.,Sergeev}
\paper{The tensor algebra of the identity representation as a module over the Lie superalgebras ${\frak Gl}(n,m)$ and Q(n)\inRus}
\journal{Mat. Sbornik}
\Vol{123}
\Year{1984}
\Pages{422-430}
\etransl{Math. USSR Sbornik}
\Vol{51}
\Year{1985}
\Pages{419-427}

\nextref Sk19
\auth{S.,Skryabin}
\paper{On the graded algebras associated with Hecke symmetries}
\quad arXiv: 1903. 06363.

\nextref Stan
\auth{R.P.,Stanley}
\book{Enumerative Combinatorics, Vol. $2$}
\publisher{Cambridge Univ. Press}
\Year{1999}

\nextref Zel
\auth{A.V.,Zelevinsky}
\book{Representations of Finite Classical Groups}
\BkSer{Lecture Notes Math.}
\BkVol{869}
\publisher{Springer}
\Year{1981}

\endreferences
\bye